\documentclass[a4paper,11pt,leqno,twoside]{article}

\usepackage{amsmath, amsthm, amssymb,bm, youngtab}
\numberwithin{equation}{section}

\makeatletter
\newif\ify@autoscale \y@autoscaletrue \def\Yautoscale#1{\ifnum #1=0
  \y@autoscalefalse\else\y@autoscaletrue\fi}
\newdimen\y@b@xdim
\newdimen\y@boxdim \y@boxdim=13pt
\def\Yboxdim#1{\y@autoscalefalse\y@boxdim=#1}
\newdimen\y@linethick    \y@linethick=.3pt
\def\Ylinethick#1{\y@linethick=#1}
\newskip\y@interspace \y@interspace=0ex plus 0.3ex
\def\Yinterspace#1{\y@interspace=#1}
\newif\ify@vcenter   \y@vcenterfalse
\def\Yvcentermath#1{\ifnum #1=0 \y@vcenterfalse\else\y@vcentertrue\fi}
\newif\ify@stdtext   \y@stdtextfalse
\def\Ystdtext#1{\ifnum #1=0 \y@stdtextfalse\else\y@stdtexttrue\fi}
\newif\ify@enable@skew   \y@enable@skewfalse
\y@vcentertrue
\def\y@vr{\vrule height0.8\y@b@xdim width\y@linethick depth 0.2\y@b@xdim}
\def\y@emptybox{\y@vr\hbox to \y@b@xdim{\hfil}}
\ify@enable@skew
 \def\y@abcbox#1{\if :#1\else
   \y@vr\hbox to \y@b@xdim{\hfil#1\hfil}\fi}
 \def\y@mathabcbox#1{\if :#1\else
   \y@vr\hbox to \y@b@xdim{\hfil$#1$\hfil}\fi}
\else
 \def\y@abcbox#1{\y@vr\hbox to \y@b@xdim{\hfil#1\hfil}}
 \def\y@mathabcbox#1{\y@vr\hbox to \y@b@xdim{\hfil$#1$\hfil}}
\fi
\def\y@setdim{%
  \ify@autoscale%
   \ifvoid1\else\typeout{Package youngtab: box1 not free! Expect an
     error!}\fi%
   \setbox1=\hbox{A}\y@b@xdim=1.6\ht1 \setbox1=\hbox{}\box1%
  \else\y@b@xdim=\y@boxdim \advance\y@b@xdim by -2\y@linethick
  \fi}
\newcount\y@counter
\newif\ify@islastarg
\def\y@lastargtest#1,#2 {\if\space #2 \y@islastargtrue
  \else\y@islastargfalse\fi}
\def\y@emptyboxes#1{\y@counter=#1\loop\ifnum\y@counter>0
  \advance\y@counter by -1 \y@emptybox\repeat}
\def\y@nelineemptyboxes#1{%
  \vbox{%
    \hrule height\y@linethick%
    \hbox{\y@emptyboxes{#1}\y@vr}
    \hrule height\y@linethick}\vspace{-\y@linethick}}
\def\yng(#1){%
  \y@setdim%
  \hspace{\y@interspace}%
  \ifmmode\ify@vcenter\vcenter\fi\fi{%
  \y@lastargtest#1,
  \vbox{\offinterlineskip
    \ify@islastarg
     \y@nelineemptyboxes{#1}
    \else
     \y@ungempty(#1)
    \fi}}\hspace{\y@interspace}}
\def\y@ungempty(#1,#2){%
  \y@nelineemptyboxes{#1}
  \y@lastargtest#2,
  \ify@islastarg
   \y@nelineemptyboxes{#2}
  \else
   \y@ungempty(#2)
  \fi}
\def\y@nelettertest#1#2. {\if\space #2 \y@islastargtrue
  \else\y@islastargfalse\fi}
\def\y@abcboxes#1#2.{%
  \ify@stdtext\y@abcbox#1\else\y@mathabcbox#1\fi%
  \y@nelettertest #2.
  \ify@islastarg\unskip%
   \ify@stdtext\y@abcbox{#2}\else\y@mathabcbox{#2}\fi%
  \else\y@abcboxes#2.\fi}
\ify@enable@skew
 \newdimen\y@full@b@xdim
 \newcount\y@m@veright@cnt
 \def\y@get@m@veright@cnt#1#2.{%
   \if :#1 \advance\y@m@veright@cnt by 1\y@get@m@veright@cnt#2.\fi}
 \let\y@setdim@=\y@setdim
 \def\y@setdim{%
   \y@setdim@ \y@full@b@xdim=\y@b@xdim
   \advance\y@full@b@xdim by 1\y@linethick}
 \def\y@m@veright@ifskew#1{
   \y@m@veright@cnt=0 \y@get@m@veright@cnt#1.
   \moveright \y@m@veright@cnt\y@full@b@xdim}
\else
 \def\y@m@veright@ifskew#1{}
\fi
\def\y@nelineabcboxes#1{%
  \y@nelettertest #1.
  \ify@islastarg
   \y@m@veright@ifskew{#1}
    \vbox{
      \hrule height\y@linethick%
      \hbox{\ify@stdtext\y@abcbox#1\else\y@mathabcbox#1\fi\y@vr}
      \hrule height\y@linethick}\vspace{-\y@linethick}
  \else
   \y@m@veright@ifskew{#1}
    \vbox{
      \hrule height\y@linethick%
      \hbox{\y@abcboxes #1.\y@vr}%
      \hrule height\y@linethick}\vspace{-\y@linethick}
  \fi}
\def\young(#1){%
  \y@setdim%
  \hspace{\y@interspace}%
  \y@lastargtest#1,
  \ifmmode\ify@vcenter\vcenter\fi\fi{%
  \vbox{\offinterlineskip
    \ify@islastarg\y@nelineabcboxes{#1}%
    \else\y@ungabc(#1)%
    \fi}}\hspace{\y@interspace}}
\def\y@ungabc(#1,#2){%
  \y@nelineabcboxes{#1}%
  \y@lastargtest#2,
  \ify@islastarg\y@nelineabcboxes{#2}%
  \else\y@ungabc(#2)%
  \fi}
\makeatother
\newcommand{\bC}{\mathbb{C}}

\newcommand{\bZ}{\mathbb{Z}}

\newcommand{\mf}[1]{\mathfrak{#1}}
\newcommand{\mr}[1]{\mathrm{#1}}
\newcommand{\mcal}[1]{\mathcal{#1}}

\newcommand{\mbb}[1]{\mathbb{#1}}

\def\PP{\mathbb{P}}

\def\sgn{\mathrm{sgn}}

\newcommand{\AV}{\mathcal{A}}

\def\({ \left( }
\def\){ \right)}

\def\Lra{\Leftrightarrow}
\def\lra{\leftrightarrow}

\theoremstyle{plain}
\newtheorem{thm}{Theorem}[section]
\newtheorem{prop}[thm]{Proposition}
\newtheorem{lem}[thm]{Lemma}
\newtheorem{cor}[thm]{Corollary}
\theoremstyle{definition}

\newtheorem{example}{Example}[section]
\newtheorem{remark}{Remark}[section]
\theoremstyle{conjecture}
\newtheorem{conj}[thm]{Conjecture}
\theoremstyle{problem}

\date{\empty}

\title{\bfseries Jucys-Murphy elements, orthogonal matrix integrals, and
Jack measures}
\author{\textsc{Sho MATSUMOTO} \\
 {\it Graduate School of Mathematics, Nagoya University.} \\
{\it  \small Furocho, Chikusa-ku, Nagoya, 
464-8602, JAPAN.} \\
\texttt{sho-matsumoto@math.nagoya-u.ac.jp}
}

\begin{document}
%>>>>>>>>>>>>>>>>>>>>>>>>>>>>>>>>>>>>>>>>>
%\setlength{\baselineskip}{15pt}

\maketitle

\begin{abstract}
We study symmetric polynomials whose variables are 
odd-numbered Jucys-Murphy elements.
They define elements of the Hecke algebra associated to
the Gelfand pair of the symmetric group with the hyperoctahedral group.
We evaluate their expansions in zonal spherical functions and
in double coset sums.
These evaluations are related to integrals of polynomial functions 
over orthogonal groups.
Furthermore, we give an extension of them, based on Jack polynomials.

\noindent
\emph{2000 Mathematics Subject Classification:} 
05E15, %Combinatorial problems concerning the classical groups
20C08, %Hecke algebras and their representations
05E05, %Symmetric functions
20C35. %Applications of group representations to physics
\end{abstract}

%\tableofcontents

\section{Introduction}

Let $S_n$ be the symmetric group.
Jucys-Murphy elements are formal sums in the group algebra $\mathbb{C}[S_{n}]$ defined by
$$
J_1=0, \qquad J_k=\sum_{i=1}^{k-1} (i \ k), \qquad k=2,3,\dots,n,
$$
where $(i \ k)$ is the transposition between $i$ and $k$.
Every two of them commute with each other.
Jucys \cite{Jucys} studied symmetric polynomials in variables $J_1,J_2,\dots,J_n$
and proved that they are central elements in $\mathbb{C}[S_n]$.
More precisely, given a symmetric function $F$, we have
$F(J_1,J_2,\dots,J_n) \in \mcal{Z}_n$, where $\mcal{Z}_n$ is the center of 
the group algebra $\mathbb{C}[S_n]$.

For a given element $F(J_1,J_2,\dots,J_n) \in \mcal{Z}_n$,
it is a natural to ask for its character expansion and class expansion.
Let $\chi^\lambda$ be the irreducible character of $S_n$ associated to 
Young diagram $\lambda$
and  $f^\lambda$ the dimension of its corresponding representation.
Jucys \cite{Jucys} obtained the character expansion 
$$
F(J_1,J_2, \dots,J_n)= \sum_{\lambda: |\lambda|=n} F(A_\lambda) \frac{f^\lambda}{n!} \chi^\lambda,
$$ 
where $|\lambda|$ is the number of boxes in the Young diagram $\lambda$
and $A_\lambda$ is the alphabet of contents $j-i$ of boxes $(i,j)$ in $\lambda$.

The computation for the class expansion of $F(J_1,J_2,\dots,J_n)$ is a more difficult problem.
Let $\mf{c}_\mu(n)$ be the sum of all permutations in $S_n$ of reduced cycle-type $\mu$.
For example, $\mf{c}_{(0)}(n)$ is the identity permutation and $\mf{c}_{(1)}(n)$
is the sum of all transpositions.
Then $\{\mf{c}_\mu(n) \ | \ |\mu|+\ell(\mu) \le n\}$ is a basis of $\mcal{Z}_n$,
where $\ell(\mu)$ is the number of rows of the diagram $\mu$.
The question is to evaluate coefficients $\AV_\mu(F,n)$ in 
$$
F(J_1,J_2,\dots,J_n) = \sum_{\mu: |\mu|+\ell(\mu) \le n} \AV_\mu(F,n) \mf{c}_\mu(n).
$$
The $\AV_\mu(F,n)$ are zero unless $\deg F \ge |\mu|$ and
polynomials in $n$. 
If $\deg F=|\mu|$,  then $\AV_\mu(F,n)$ is independent of $n$.
Lascoux and Thibon \cite{LT} studied the coefficients for power-sum symmetric functions $F=p_k$,
and Fujii et al. \cite{FKMO} expressed $\AV_{(0)}(p_k,n)$ as an explicit polynomial in
binomail coefficients $\binom{n}{m}$.
Matsumoto and Novak \cite{MN} gave a combinatorial explicit expression of 
$\AV_\mu(m_\lambda,n)$ with $|\lambda|=|\mu|$,
where $m_\lambda$ is the monomial symmetric function.

The coefficients $\AV_\mu(h_k,n)$ for complete symmetric functions $h_k$ are more interesting.
They appear in expansions of unitary matrix integrals.
Let $U(N)$ be the group of $N \times N$ unitary matrices $g=(g_{ij})_{1 \le i,j \le N}$, 
equipped with its normalized Haar measure $d g$.
Consider integrals of the form
$$
\int_{g \in U(N)} g_{i_1 j_1} g_{i_2 j_2} \cdots g_{i_n j_n} 
\overline{g_{i_1' j_1'} g_{i_2' j_2'} \cdots g_{i_n' j_n'}} dg
$$
where $i_k, i_k',j_k,j_k'$ are in $\{1,2,\dots, N\}$
and  $N \ge n$.
The Weingarten calculus for unitary groups
developed in \cite{W,C,CS} states that
those integrals are given by a sum of Weingarten functions
$$
\mr{Wg}^{U(N)}_n(\sigma)= \int_{g \in U(N)} \prod_{k=1}^n g_{kk} \overline{g_{k \sigma(k)}} d g,
\qquad \sigma \in S_{n}.
$$
In \cite{Novak} (see also \cite{MN}), 
a remarkable connection between $\mr{Wg}^{U(N)}_n$ and Jucys-Murphy elements is
discovered.
Specifically, the Weingarten function is given as a generating function of $h_k(J_1,\dots,J_n)$:
$$
\sum_{\sigma \in S_n} \mr{Wg}^{U(N)}_n(\sigma) \sigma = \sum_{k=0}^\infty (-1)^k N^{-n-k} h_k(J_1,J_2,\dots,J_n),
$$
or equivalently
$$
\mr{Wg}^{U(N)}_n(\sigma) = \sum_{k=0}^\infty (-1)^k N^{-n-k} \AV_\mu(h_k,n),
$$
where $\mu$ is the reduced cycle-type of $\sigma$.
Thus unitary matrix integrals are evaluated by observing symmetric functions in Jucys-Murphy elements.

The main purpose of the present paper is to study their analogues for orthogonal groups.
In the orthogonal group case, the elements
$F(J_1,J_3,\dots,J_{2n-1}) \cdot P_n$ are needed instead of $F(J_1,J_2,\dots,J_n)$.
Here $P_n=\sum_{\zeta \in H_n} \zeta$ is an element of $\mathbb{C}[S_{2n}]$ 
and $H_n$ is the hyperoctahedral group realized in $S_{2n}$.
We will prove first that $F(J_1,J_3,\dots,J_{2n-1}) \cdot P_n$ belongs to
the Hecke algebra $\mcal{H}_n$ associated with the Gelfand pair $(S_{2n},H_n)$.
The Hecke algebra $\mcal{H}_n$ has two kinds of natural basis 
$\{\omega^\lambda \ | \ |\lambda|=n\}$ and $\{\psi_\mu(n) \ | \ |\mu|+\ell(\mu) \le n \}$,
where the $\omega^\lambda$ are zonal spherical functions and 
the $\psi_\mu(n)$ are sums over double cosets of the form $H_n \sigma H_n$.
As in the unitary group case, 
it is natural to ask for expansions of $F(J_1,J_3,\dots,J_{2n-1}) \cdot P_n$
in zonal spherical functions $\omega^\lambda$ or in double coset sums $\psi_\mu(n)$.
We will prove that the expansion in $\omega^\lambda$ is given by 
$$
F(J_1,J_3,\dots,J_{2n-1}) \cdot P_n= 
\sum_{\lambda: |\lambda|=n} F(A_\lambda') \frac{f^{2\lambda}}{(2n-1)!} \omega^\lambda,
$$
where $A_\lambda'$ is the alphabet of modified contents $2j-i-1$.
Our main purpose is to obtain some properties for coefficients $\AV'_\mu(F,n)$ defined via
$$
F(J_1,J_3,\dots,J_{2n-1} ) \cdot P_n =
\sum_{\mu: |\mu|+\ell(\mu) \le n} \AV_\mu'(F,n) \mf{\psi}_\mu(n).
$$
In general, the $\AV_\mu'(F,n)$ are different from 
$\AV_\mu(F,n)$.
For example, $\AV_{(1)}(h_3,n)= \tfrac{1}{2}n^2+\tfrac{3}{2}n-4$ but 
$\AV'_{(1)}(h_3,n)=n^2+3n-7$.
However, by observing deep combinatorics of perfect matchings,
 we will prove that,
if $\deg F=|\mu|$, they coincide as $\AV_\mu(F,n)=\AV_\mu'(F,n)$,
and are independent of $n$.

Like in the unitary group case,
coefficients $\AV_\mu'(h_k,n)$ are involved in orthogonal matrix integrals.
Let $O(N)$ be the orthogonal group of degree $N$ and
$d g$  its normalized Haar measure. 
Then, for example, we obtain
$$
\int_{g \in O(N)} g_{11}^2 g_{22}^2 \cdots g_{nn}^2 d g =
\sum_{k=0}^\infty (-1)^{k} N^{-n-k} \AV_{(0)}'(h_k,n).
$$
Via the Weingarten calculus for orthogonal groups developed in \cite{CS,CM,ZJ},
we establish the connection between orthogonal matrix integrals and Jucys-Murphy elements.

Furthermore, we introduce an $\alpha$-extension of $\AV_\mu(F,n)$ and $\AV'_\mu(F,n)$.
Let $\alpha$ be a positive real number.
We define the value $\AV_\mu^{(\alpha)}(F,n)$ as an average with respect to 
the Jack measure. 
The Jack measure  is a probability measure on Young diagrams
and is a deformation of the Plancherel measure for symmetric groups.
Its definition is based on Jack polynomial theory and 
the connections between them and random  matrix theory are much studied, see \cite{Mat, BO} and their references.
From symmetric function theory, we can see $\AV_\mu(F,n)=\AV^{(1)}_\mu(F,n)$
and $\AV_\mu'(F,n)= \AV^{(2)}_\mu(F,n)$ for any symmetric function $F$ and partition $\mu$.
Also $\AV_\mu^{(1/2)}(F,n)$ are important and related to a twisted Gelfand pair.

We evaluate $\AV_\mu^{(\alpha)}(e_k,n)$ for elementary symmetric functions $e_k$.
Also, by applying shifted symmetric function theory developed in 
\cite{KOO,LassalleSomeIdentities, LassalleCumulant, Olshanski},
we prove that
the $\AV^{(\alpha)}_\mu(F,n)$ are polynomials in $n$.
We could not obtain any strong results for $\AV_\mu^{(\alpha)}(F,n)$. 
Our appoarch is experimental but the author believes that it is fascinating and 
applicable in futurer research.

The present paper is constructed as follows:
In Section 2 we review necessary notations and fundamental properties. 
In Section 3, we evaluate $e_k(J_1,J_3,\dots,J_{2n-1}) \cdot P_n$ explicitly and
prove that $F(J_1,J_3,\dots,J_{2n-1}) \cdot P_n$ is an element of the Hecke algebra $\mcal{H}_n$.
In Section 4, we give the expansion of $F(J_1,J_3,\dots,J_{2n-1}) \cdot P_n$
in zonal spherical functions.
In Section 5 an 6, we give some properties of $\AV_\mu'(F,n)$.
Specifically, we prove that $\AV_\mu'(F,n)$ coincides with
$\AV_\mu(F,n)$ if $|\mu|= \deg F$.
As we mentioned, such an equality does not hold for  $|\mu|< \deg F$.
In Section 7, we see the connection to orthogonal matrix integrals.
In Section 8, we study Jack's $\alpha$-deformations $\AV_\mu^{(\alpha)}(F,n)$.
In the final section, Section 9,
we give some examples and suggest four conjectures.

\begin{remark}
Since a primary version of this paper was released,
all of our four conjectures given in Subsection \ref{subsec:Open}
have been actively studied by some other reseachers.
We would be able to see their proofs very soon.
\end{remark}

\section{Preparations} \label{Sec:Preparations}

We use the notations of Macdonald.
See Chapter I and VII in his book \cite{Mac}.

\subsection{Partitions and contents} \label{subsec:partitions}

A {\it partition} $\lambda=(\lambda_1,\lambda_2,\dots)$ is a weakly decreasing sequence of 
nonnegative integers such that 
its {\it length} $\ell(\lambda):=|\{ i \ge 1 \ | \ \lambda_i >0\}|$ is finite.
We write the {\it size} of $\lambda$ by $|\lambda|:= \sum_{i \ge 1} \lambda_i$.
If $|\lambda|=n$, we say $\lambda$ to be  a partition of $n$ and write $\lambda \vdash n$.

We often identify $\lambda$ with its {\it Young diagram}
$Y(\lambda):=\{\square =(i,j) \in \bZ^2 \ | \ 1 \le i \le \ell(\lambda), \ 1 \le j \le \lambda_i \}$.
If $\square =(i,j) \in Y(\lambda)$, we say that $\square$ is a {\it box} of $\lambda$ and write
$\square \in \lambda$ shortly.
The {\it content} of $\square=(i,j) \in \lambda$ is defined by $c(\square):=j-i$.
Also we use its analogy $c'(\square):=2j-i-1$.
Let $A_\lambda$ and $A_\lambda'$ be the alphabet with $|\lambda|$ elements given by
$$
A_\lambda=\{c(\square) \ | \ \square \in \lambda\}, \qquad
A_\lambda'=\{c'(\square) \ | \ \square \in \lambda\}.
$$
For example, $A_{(2,2,1)}=\{1,0,0,-1,-2\}$ and $A_{(2,2,1)}'=\{2,1,0,-1,-2\}$.

For each $i \ge 1$, we write the multiplicity of $i$ in $\lambda$ by 
$m_i(\lambda)=|\{j \ge 1 \ | \ \lambda_j=i\}|$.
We sometimes write $\lambda$ as 
$(\dots,3^{m_3(\lambda)}, 2^{m_2(\lambda)},1^{m_1(\lambda)})$.
For example, $\lambda=(2,1,1,1)= (2,1^3)$.
Define
\begin{equation} \label{eq:zlambda}
z_\lambda= \prod_{i \ge 1} i^{m_i(\lambda)} m_i(\lambda)!.
\end{equation}

Let $\lambda,\mu$ be partitions.
We define $\lambda+\mu$ to be the sequence of $\lambda_i+\mu_i$:
$(\lambda+\mu)_i=\lambda_i+\mu_i$.
Also we define $\lambda \cup \mu$ to be the partition whose parts are those of $\lambda$
and $\mu$, arranged in decreasing order.
In general, given partitions $\lambda^{(1)},\lambda^{(2)},\dots,\lambda^{(k)}$,
we define $\lambda^{(1)} \cup \lambda^{(2)} \cup \dots \cup \lambda^{(k)}$ in the same way.

For a partition $\lambda$ with $\ell(\lambda)=l$, 
we define its {\it reduction} $\tilde{\lambda}$ by 
$\tilde{\lambda}=(\lambda_1-1, \lambda_2-1,\dots,\lambda_{l}-1)$.
For each $n \ge 1$, the map $\lambda \mapsto \tilde{\lambda}$ gives a bijection
from the set
$\{\lambda \ | \ |\lambda|= n\}$ to 
$\{ \mu \ | \ |\mu|+\ell(\mu) \le n \}$.
Indeed, its inverse map is given by $\mu 
\mapsto \mu+(1^{n-|\mu|})=:\lambda$.
Then $|\lambda|-\ell(\lambda)=|\mu|$.

\subsection{Symmetric functions}
 \label{subsec:SymmetricFunctions}

Let $x=(x_1,x_2,\dots)$ be an infinite sequence of indeterminates,
and $\mbb{S}$ the algebra of symmetric functions with complex coefficients in variables $x$.

Given a partition $\lambda$, the \emph{monomial symmetric function} $m_\lambda$ 
is defined by
$$
m_\lambda(x)= \sum_{(\alpha_1,\alpha_2,\dots)} x_1^{\alpha_1}
x_2^{\alpha_2}\cdots ,
$$
summed over all distinct permutations $(\alpha_1,\alpha_2,\dots)$ of
$(\lambda_1,\lambda_2,\dots)$.
Denote by $e_k, p_k$, and $h_k$ the 
\emph{elementary}, \emph{power-sum}, and \emph{complete symmetric functions},
respectively.
Namely, 
\begin{align*}
e_k(x)=& m_{(1^k)}(x)= \sum_{i_1 < i_2 < \cdots <i_k} x_{i_1} x_{i_2} \cdots x_{i_k}, \\
p_k(x)=&m_{(k)}(x)= x_1^k+x_2^k+\cdots, \\
h_k(x)=& \sum_{\lambda \vdash k} m_\lambda(x)= \sum_{ i_1 \le i_2 \le \cdots \le i_k}x_{i_1} x_{i_2} \cdots x_{i_k}.
\end{align*}
Also we put 
$e_{\lambda}=\prod_{i=1}^{\ell(\lambda)} e_{\lambda_i}$, 
and similarly for $p_\lambda$ and $h_\lambda$.
For convenience, we set $m_{(0)}=e_0=h_0=1$.

For finite variables $(x_1,x_2,\dots,x_n)$, 
the monomial symmetric function (or polynomial) $m_\lambda(x_1,x_2,\dots)$ is zero
unless $\ell(\lambda) \le n$.

The degree of $m_\lambda$ is naturally defined to be $\deg m_\lambda= |\lambda|$.

The fundamental theorem for symmetric functions says that
any symmetric function $F$ is given by a polynomial in 
$e_1,e_2,\dots$ and that the $e_k$ are algebraically independent.
We can replace $e_k$ by $p_k$ in this statement.

\subsection{Hyperoctahedral groups} \label{subsec:Hyperoctahedral}

We recall a Gelfand pair $(S_{2n},H_n)$. The detail is seen in \cite[VII.2]{Mac}.

Let $S_{n}$ be the symmetric group on $\{1,2,\dots, n\}$. 
Let $\bC[S_{n}]$ denote the algebra of all complex-valued functions $f$
on $S_n$ under convolution 
$(f_1 \cdot f_2)(\sigma)= \sum_{\tau \in S_{n}} f_1(\sigma \tau^{-1}) f_2(\tau)$.
This is identified with the algebra of formal $\mathbb{C}$-linear sums of 
permutations with multiplication 
$(\sum_{\sigma} f_1(\sigma) \sigma)(\sum_{\tau} f_2(\tau) \tau) =
\sum_{\pi} \(\sum_{\sigma} f_1(\sigma) f_2(\sigma^{-1} \pi)  \) \pi$.

A permutation $\sigma$ in $S_n$ is regarded as a permutation in $S_{n+1}$ fixing 
the letter $n+1$. Thus $\bC[S_{n}] \subset  \bC[S_{n+1}]$.

Let $H_n$ be the \emph{hyperoctahedral group}, which 
is a subgroup of $S_{2n}$ generated by
transpositions $(2i-1 \ 2i)$, $(1 \le i \le n)$, and by
double transpositions $(2i-1 \ 2j-1) (2i \ 2j)$, $(1 \le i<j \le n)$.
Equivalently, $H_n$ is the centralizer of $(1 \ 2)(3 \ 4) \cdots (2n-1 \ 2n)$ in $S_{2n}$.
Then the pair $(S_{2n},H_n)$ is a Gelfand pair.
Let $P_n$ be the sum of all elements of $H_n$ in $\mathbb{C}[S_{2n}]$:
$$
P_n=\sum_{\zeta \in H_n}\zeta. 
$$

Consider the double cosets $H_n \sigma H_n$ in $S_{2n}$.
These cosets are indexed by partitions of $n$, that is,
\begin{equation} \label{eq:bothsidedecomposition}
S_{2n} = \bigsqcup_{\rho \vdash n} H_{\rho}, 
\end{equation}
where each $H_\rho$ is a double coset.
The permutation $\sigma \in S_{2n}$ is said to be of {\it coset-type} $\rho$
and written as $\Xi_n(\sigma)= \rho$
if $\sigma \in H_\rho$.

Also the partition $\Xi_n(\sigma)$ is graphically defined as follows.
Consider the graph $\Gamma(\sigma)$ whose vertex set is $\{1,2,\dots,2n\}$ and whose
edge set consists of $\{2i-1,2i\}$ and $\{\sigma(2i-1),\sigma(2i)\}$, $1 \le i \le n$.
We think of the edges $\{\sigma(2i-1),\sigma(2i)\}$ as blue, and the others as red.
Then $\Gamma(\sigma)$ has some connected components of even lengths 
$2 \rho_1 \ge 2 \rho_2 \ge \cdots$.
Thus $\sigma$ determines a partition $\rho=(\rho_1,\rho_2,\dots)$ of $n$.
The $\rho$ is nothing but the coset-type $\Xi_n(\sigma)$.

Two permutations $\sigma_1,\sigma_2 \in S_{2n}$ have the same coset-type if
and only if $H_n \sigma_1 H_n= H_n \sigma_2 H_n$.
The set $H_\rho$ consists of permutations in $S_{2n}$ of coset-type $\rho$.
Given $\sigma \in S_{2n}$, we let $\nu_n(\sigma)$ to be the length of 
the partition $\Xi_n(\sigma)$: $\nu_n(\sigma) = \ell(\Xi_n(\sigma))$.

\begin{example} \label{ex:GraphGamma}
For $\sigma=\(\begin{smallmatrix} 1 & 2 & 3 & 4 & 5 & 6 & 7 & 8 & 9 & 10 \\
5 & 1 & 4 & 10 & 3 & 9 & 7 & 6 & 2 & 8 \end{smallmatrix} \) \in S_{10}$,
its graph $\Gamma(\sigma)$ has two connected components $\Gamma^{(1)}$ and $\Gamma^{(2)}$:
$$
\Gamma^{(1)}: 1 \Longleftrightarrow 5  \longleftrightarrow 6
\Longleftrightarrow 7 \longleftrightarrow 8 \Longleftrightarrow 2
\longleftrightarrow 1, \qquad
\Gamma^{(2)}: 3 \Longleftrightarrow 9 \longleftrightarrow 10 \Longleftrightarrow 4 
\longleftrightarrow 3.
$$
Here ``$i \Longleftrightarrow j$'' means that
a blue edge connects 
the $i$-th vertex with the $j$-th vertex,
whereas ``$p \longleftrightarrow q$'' means that
a red edge connects 
the $p$-th vertex with the $q$-th vertex.
Equivalently, there exists an integer $k$ such that
$\{i,j\}=\{\sigma(2k-1), \sigma(2k)\}$ (resp. $\{p,q\}=\{2k-1,2k\}$).
In this example, the component $\Gamma^{(1)}$ 
and $\Gamma^{(2)}$ has $6$ and $4$ vertices, respectively,
and hence $\Xi_5(\sigma)=(3,2)$ and $\nu_5(\sigma)=2$.
\end{example}

\subsection{Perfect matchings} \label{subsec:PM}

Let $\mcal{M}(2n)$ be the set of all \emph{perfect matchings} on $\{1,2,\dots,2n\}$.
Each perfect matching $\mf{m}$ in $\mcal{M}(2n)$ is uniquely expressed by the form
\begin{equation} \label{eq:ExpressionPairPartition}
\left\{ \{\mf{m}(1),\mf{m}(2) \}, \{\mf{m}(3),\mf{m}(4) \}, \dots, \{\mf{m}(2n-1),\mf{m}(2n) \}
\right\}
\end{equation}
with $\mf{m}(2i-1) < \mf{m}(2i)$ for $1 \le i \le n$ and with
$\mf{m}(1) <\mf{m}(3) < \cdots < \mf{m}(2n-1)$.
We call each $\{\mf{m}(2i-1),\mf{m}(2i)\}$ a \emph{block} of $\mf{m}$.
A block  of the form $\{2i-1,2i\}$ is said to be \emph{trivial}. 
We embed the set $\mcal{M}(2n)$ into $S_{2n}$ via the mapping
\begin{equation} \label{eq:identificationPM}
\mcal{M}(2n) \ni \mf{m} \mapsto
\begin{pmatrix}
1 & 2 & 3 & 4 & \cdots & 2n \\
\mf{m}(1) & \mf{m}(2) & \mf{m}(3) & \mf{m}(4) & \cdots & \mf{m}(2n) 
\end{pmatrix} \in S_{2n}.
\end{equation}
In particular, the graph $\Gamma(\mf{m})$, the coset-type $\Xi_n(\mf{m})$, 
and the value $\nu_n(\mf{m})$ are defined.
Note that  $\Gamma(\mf{m})=\Gamma(\mf{n})$ if and only if
$\mf{m}=\mf{n}$.

A perfect matching $\mf{n}$ in $\mcal{M}(2n-2)$ is regarded as an element of $\mcal{M}(2n)$ by
adding the trivial block $\{2n-1,2n\}$:
$$
\mcal{M}(2n-2) \ni \mf{n} \mapsto \mf{n} \sqcup \{\{2n-1,2n\}\} \in \mcal{M}(2n).
$$
Thus we think as  $\mcal{M}(2n-2) \subset \mcal{M}(2n)$.

It is well known that 
$\mcal{M}(2n)$ is the set of all representatives of the right cosets $\sigma H_n$
of $H_n$ in $S_{2n}$: 
\begin{equation} \label{eq:rightcosets}
S_{2n} = \bigsqcup_{\mf{m} \in \mcal{M}(2n)} \mf{m} H_n.
\end{equation}

\subsection{Characters and zonal spherical functions} \label{subsec:CharZonal}

Given a partition $\lambda \vdash n$, we denote by $\chi^\lambda$ the irreducible character
of $S_{n}$. The set $\{\chi^\lambda \ | \ \lambda \vdash n \}$ is a basis of the center of the group algebra $\mathbb{C}[S_{n}]$.
Let $\mr{id}_n$ denote the identity permutation in $S_n$ and let
$f^\lambda:= \chi^\lambda(\mr{id}_n)$.
Thus the number $f^\lambda$ is the dimension of the irreducible representation
of $S_{n}$ with character $\chi^\lambda$.
Equivalently, $f^\lambda$ is the number of standard Young tableaux of shape $\lambda$, 
see e.g. \cite{Sagan}.

For each partition $\lambda$ of $n$, we define the {\it zonal spherical function} of the 
Gelfand pair $(S_{2n},H_n)$ by
\begin{equation}
\omega^\lambda (\sigma) = \frac{1}{2^n n!} \sum_{\zeta \in H_n} \chi^{2\lambda}(\sigma \zeta),
\qquad \sigma \in S_{2n},
\end{equation}
where $2\lambda=\lambda+\lambda=(2\lambda_1,2\lambda_2,\dots)$.
If we regard $\omega^\lambda$ as an element of $\bC[S_{2n}]$,
we can express $\omega^\lambda= \frac{1}{2^n n!} \chi^{2\lambda} P_n=
\frac{1}{2^n n!} P_n \chi^{2\lambda}$.
These functions are $H_n$-biinvariant functions on $S_{2n}$ and 
 constant on each double coset $H_\rho$.
Denote by $\omega^\lambda_\rho$ the value of $\omega^\lambda$ at $H_\rho$.

Let $\mcal{H}_n$ be the Hecke algebra associated with the Gelfand pair $(S_{2n},H_n)$:
$$
\mcal{H}_n =\{ f : S_{2n} \to \mathbb{C} \ | \ \text{$f$ is constant on each $H_\rho$ ($\rho \vdash n$)} \}.
$$
Since $(S_{2n},H_n)$ is a Gelfand pair,
this algebra is commutative with respect to the convolution product.
We often regard $\mcal{H}_n$ as a subspace of $\mathbb{C}[S_{2n}]$.
There are two natural bases of $\mcal{H}_n$;
one is $\{\omega^\lambda \ | \ \lambda \vdash n\}$ and another is
$\{\phi_\rho \ | \ \rho \vdash n\}$.
Here the $\phi_\rho$ are double-coset sum functions
$$
\phi_\rho= \sum_{\sigma \in H_\rho} \sigma.
$$ 
Note that $\phi_{(1^n)}= P_n$.

\section{Analogue of Jucys' result}

Define the Jucys-Murphy elements $J_k$.
They are commuting elements in  group algebras of symmetric groups, given by
$J_1=0$ and by
$$
J_k=(1 \ k) +(2 \ k)+ \cdots +(k-1 \ k) \qquad \text{for $k = 2,3,\dots$}.
$$

Jucys \cite{Jucys} obtains an exact expression for 
$e_k(J_1,J_2,\dots,J_n)$,
where $e_k$ is the elementary symmetric function.
His result is the following identity:
$$
e_k(J_1,J_2,\dots,J_n)= \sum_{\pi} \pi
$$
summed over all permutations $\pi$ in $S_{n}$ with exactly $n-k$ cycles (including trivial cycles).
The following proposition is an analogue of this identity, and was
essentially obtained by Zinn-Justin \cite{ZJ} very recently.
Our proof is an analogue of Jucys' proof in \cite{Jucys}.

\begin{prop} \label{prop-Jucys}
For any $k$ and $n$, we have
\begin{equation} \label{eq:ElementaryJM1}
e_k(J_1,J_3,\dots,J_{2n-1}) \cdot  P_n = 
\sum_{\begin{subarray}{c} \mf{m} \in \mcal{M}(2n) \\
\nu_n(\mf{m})= n-k \end{subarray}} \mf{m} P_n.
\end{equation}
\end{prop}

\begin{proof}
First observe that \eqref{eq:ElementaryJM1} holds true when $k=0$
because $\mf{m}=\{\{1,2\},\{3,4\},\dots,\{2n-1,2n\}\}$ is 
the unique element in $\mcal{M}(2n)$ satisfying $\nu_n(\mf{m})=n$.
Also, when $k \ge n$, both sides are zero.

We proceed by induction on $n$.
Let $n >1$ and suppose that the claim holds for 
$e_k(J_1,J_3,\dots,J_{2n-3}) \cdot P_{n-1}$ with any $k \ge 0$.
Using identities 
$e_{k}(x_1,x_2,\dots,x_n)= e_{k}(x_1,\dots,x_{n-1})
+ x_n e_{k-1}(x_1,\dots,x_{n-1})$ and 
$P_{n-1} P_n =|H_{n-1}| P_n$,
we see that
\begin{align*}
& e_k(J_1,J_3,\dots,J_{2n-1})  P_{n} \\
=& e_k (J_1,J_3,\dots,J_{2n-3})  P_n + J_{2n-1} e_{k-1}(J_1,J_3,\dots,J_{2n-3})P_n \\
=& \frac{1}{|H_{n-1}|} e_{k}(J_1,J_3,\dots,J_{2n-3}) P_{n-1} P_n 
+ \frac{1}{|H_{n-1}|} J_{2n-1} e_{k-1} (J_1,J_3,\dots,J_{2n-3}) P_{n-1} P_n.
\end{align*}
The induction assumption gives
\begin{equation} \label{eq:InductionElementary}
e_k(J_1,J_3,\dots,J_{2n-1})  P_{n}=
\sum_{\begin{subarray}{c} \mf{n} \in \mcal{M}(2n-2) \\
\nu_{n-1} (\mf{n})= n-1-k \end{subarray}} \mf{n} P_n 
+ \sum_{t=1}^{2n-2}  \sum_{\begin{subarray}{c} \mf{n} \in \mcal{M}(2n-2) \\
\nu_{n-1}(\mf{n})= n-k \end{subarray}} (t \ 2n-1) \mf{n}P_n.
\end{equation}

Recall the natural inclusion 
$\mcal{M}(2n-2) \ni \mf{n} \mapsto \mf{n} \sqcup \{\{2n-1, 2n\}\} \in \mcal{M}(2n)$.
We have $\nu_n(\mf{n} \sqcup \{\{2n-1, 2n\}\} )=
\nu_{n-1} (\mf{n})+1$,
and hence the first summand on the right hand side of \eqref{eq:InductionElementary}
coincides with $\sum_{\mf{m}} \mf{m} P_n$,
summed over $\mf{m} \in \mcal{M}(2n)$ having the block $\{2n-1,2n\}$ with 
$\nu_{n}(\mf{m})=n-k$.

Next, let us see  $ (t \ 2n-1) \mf{n}P_n$ for $\mf{n} \in \mcal{M}(2n-2)$ and
$1 \le t \le 2n-2$.
Denote by $t_{\mf{n}}$ the index in $\{1,2,\dots,2n-2\}$,
determined by $\{t_{\mf{n}},t\} \in \mf{n}$. 
We define an element $\mf{n}_t$ in $\mcal{M}(2n)$
by removing $\{t_{\mf{n}},t\}$ from $\mf{n}$ and by
adding two new blocks $\{t_{\mf{n}},2n-1\}$ and $\{t,2n\}$:
$$
\mf{n}_t = (\mf{n} \setminus \{\{t_{\mf{n}},t\}\}) \cup \{ \{ t_{\mf{n}},2n-1\},
\{t,2n\} \}.
$$
Then it is easy to see $\mf{n}_t H_n =
(t \ 2n-1)\mf{n}H_n$, and therefore 
$$
\mf{n}_t P_n= (t \ 2n-1) \mf{n}  P_n.
$$
Moreover, we can see 
$$
\nu_{n-1}(\mf{n})= \nu_n(\mf{n}_t).
$$
In fact, consider graphs $\Gamma(\mf{n}), \Gamma(\mf{n}_t)$
defined in Subsection \ref{subsec:Hyperoctahedral}.
We use the notation of Example \ref{ex:GraphGamma}.
The graph $\Gamma(\mf{n}_t)$ can be obtained from $\Gamma(\mf{n})$ if we replace 
an edge $t_{\mf{n}} \Longleftrightarrow t$  in $\Gamma(\mf{n})$ 
by the path $t_{\mf{n}} \Longleftrightarrow 2n-1 \longleftrightarrow
2n \Longleftrightarrow t$.
This means that the numbers of components of $\Gamma(\mf{n})$ and $\Gamma(\mf{n}_t)$
coincide, i.e., $\nu_{n-1}(\mf{n})= \nu_n(\mf{n}_t)$.

The observation in the previous paragraph implies that
for each $t$
the sum  $\sum_{\mf{n}} (t \ 2n-1) \mf{n}P_n$
on \eqref{eq:InductionElementary}
coincides with $\sum_{\mf{m}} \mf{m} P_n$,
summed over $\mf{m} \in \mcal{M}(2n)$ having the block $\{t,2n\}$ with 
$\nu_{n}(\mf{m})=n-k$.

It follows that the expression \eqref{eq:InductionElementary} 
is $\sum_{\mf{m}} \mf{m} P_n$,
summed over all $\mf{m} \in \mcal{M}(2n)$ with 
$\nu_{n}(\mf{m})=n-k$.
Thus \eqref{eq:ElementaryJM1} is proved.
\end{proof}

Recall that 
$H_\rho$ is the double coset $H_n \sigma H_n$ of permutations
of coset-type $\rho$, and that $\phi_\rho$
is the formal sum over $H_\rho$ in $\mathbb{C}[S_{2n}]$.

\begin{cor} \label{cor1-Jucys}
For each $0 \le k <n$, we have
\begin{equation} \label{eq:ElementaryJM2}
e_k(J_1,J_3,\dots,J_{2n-1}) \cdot P_n 
= \sum_{\begin{subarray}{c} \rho \vdash n \\
\ell(\rho)= n-k \end{subarray}} \phi_\rho.
\end{equation}
This belongs to $\mcal{H}_n$.
Moreover, 
$e_k(J_1,J_3,\dots,J_{2n-1}) \cdot P_n 
= P_n \cdot e_k(J_1,J_3,\dots,J_{2n-1})$.
\end{cor}

\begin{proof}
By decompositions \eqref{eq:rightcosets} and \eqref{eq:bothsidedecomposition},
the right hand side on \eqref{eq:ElementaryJM1} equals
$$
\sum_{\begin{subarray}{c} \mf{m} \in \mcal{M}(2n) \\
\nu_n(\mf{m})= n-k \end{subarray}} \sum_{\zeta \in H_n} \mf{m} \zeta
=\sum_{\begin{subarray}{c} \sigma \in S_{2n} \\
\nu_n(\sigma)= n-k \end{subarray}} \sigma 
=\sum_{\begin{subarray}{c} \rho \vdash n \\
\ell(\rho)= n-k \end{subarray}} \sum_{\sigma \in H_\rho} \sigma,
$$
which implies \eqref{eq:ElementaryJM2}. Let $\iota: \bC[S_{2n}] \to \bC[S_{2n}]$
be the linear extension of the anti-isomorphism $S_{2n} \ni\sigma \mapsto \sigma^{-1} \in S_{2n}$.
By \eqref{eq:ElementaryJM2}
and by the fact that $H_n \sigma H_n=H_n \sigma^{-1} H_n$ for any $\sigma \in S_{2n}$, 
the element $e_k(J_1,J_3,\dots,J_{2n-1}) \cdot P_n$
is invariant under $\iota$. 
However, $\iota (e_k(J_1,J_3,\dots,J_{2n-1}) \cdot P_n) =\iota (P_n) \cdot \iota(e_k(J_1,J_3,\dots,J_{2n-1})) = P_n \cdot e_k(J_1,J_3,\dots,J_{2n-1})$.
\end{proof}

Now the fundamental theorem on symmetric polynomials  gives

\begin{cor} \label{cor:SymPolyHecke}
For any symmetric function $F$ and positive integer $n$, 
$$
F(J_1,J_3,\dots,J_{2n-1}) \cdot P_n = P_n \cdot F(J_1,J_3,\dots,J_{2n-1}),
$$
which belongs to $\mcal{H}_n$.
\end{cor}

We are interested in the expansion of $F(J_1,J_3,\dots,J_{2n-1}) \cdot P_n$
with respect to basis $\omega^\lambda$'s or $\phi_\rho$'s of $\mcal{H}_n$.

\section{Spherical expansion}

Our purpose in this section is to obtain the expansion of $F(J_1,J_3,\dots,J_{2n-1}) \cdot P_n$ 
in zonal spherical functions $\omega^\lambda$.

Given $F \in \mbb{S}$ and $\lambda \vdash n$, we put
$$
F(A_\lambda') := F(x_1,x_2,\dots,x_n,0,0,\dots)|_{\{x_1,x_2,\dots,x_n\} =A_\lambda'}
$$ 
where $A_\lambda'$ was defined in Subsection \ref{subsec:partitions}.

\begin{thm} \label{thm:fOmegaContent}
For any $\lambda \vdash n$ and symmetric function $F$, 
$$
F(J_1,J_3,\dots,J_{2n-1}) \cdot \omega^\lambda = 
 \omega^\lambda  \cdot F(J_1,J_3,\dots,J_{2n-1}) = F(A_\lambda') \omega^\lambda.
$$
\end{thm}

\begin{proof}
In this proof, we suppose readers are familiar with standard Young tableaux, see e.g. \cite{Sagan}.
For a partition $\mu$, denote by $\mr{SYT}(\mu)$
the set of all standard Young tableaux of shape $\mu$.
For each standard Young tableau $T \in \mr{SYT}(\mu)$, 
let $e_T \in \bC[S_{|\mu|}]$ be Young's orthogonal idempotent.
Their definition and properties are seen in \cite{Garsia}.
We use well-known identities
$$
J_k \cdot e_T= e_T \cdot J_k=c(T_k) e_T, \qquad 
 \sum_{T \in \mr{SYT}(\mu)} e_T =  \frac{f^\mu}{|\mu|!} \chi^\mu.
$$
Here $T_k$ is the box $\square=(i_k,j_k)$ in $T$ labelled by $k$
and $c(T_k)$ is its content $j_k-i_k$.
Note $f^\mu = |\mr{STY}(\mu)|$.
These identities imply that, for each $\mu \vdash 2n$,
\begin{equation} \label{eq:YngIdenF}
F(J_1,J_3,\dots,J_{2n-1}) \cdot \chi^\mu=
\frac{(2n)!}{f^\mu} \sum_{T \in \mr{SYT}(\mu)} F(c(T_1),c(T_3),\dots,c(T_{2n-1}))e_T.
\end{equation}

Let $\lambda \vdash n$. 
Given $S =(S[i,j])_{(i,j) \in \lambda} \in \mr{SYT}(\lambda)$, we define the standard Young tableau $S'=(S'[i,j])_{(i,j) \in 2\lambda} \in \mr{SYT}(2\lambda)$ by
$$
S'[i,2j-1]= 2S[i,j]-1, \qquad S'[i,2j]= 2S[i,j], \qquad (i,j) \in \lambda.
$$
Here $S[i,j]$ stands for the entry of $S$ in the box $(i,j)$.
For example, given $S={\footnotesize \young(13,24)} \in \mr{SYT}((2,2))$, 
we have $S'={\footnotesize \young(1256,3478)} \in \mr{SYT}((4,4))$.
Proposition 4 in \cite{ZJ} claims that,
given a standard Young tableau $T$ with $2n$ boxes,
$P_n e_T$ is zero unless there is a standard tableau $S$ with $n$ boxes such that
$T=S'$.
We have $c(S'_{2k-1})=c'(S_k)$
by the construction of $S'$.
Hence it follows by \eqref{eq:YngIdenF} that 
\begin{align*}
&F(J_1,J_3,\dots,J_{2n-1}) \cdot P_n \cdot \chi^{2\lambda}
= P_n \cdot F(J_1,J_3,\dots,J_{2n-1})  \cdot \chi^{2\lambda} \\
=& \frac{(2n)!}{f^{2\lambda}} \sum_{T \in \mr{SYT}(2\lambda)} 
F(c(T_1),c(T_3),\dots,c(T_{2n-1})) P_n \cdot e_T \\
=& \frac{(2n)!}{f^{2\lambda}} \sum_{S \in \mr{SYT}(\lambda)} 
F(c'(S_1),c'(S_2),\dots,c'(S_n)) P_n \cdot e_{S'} \\
=& \frac{(2n)!}{f^{2\lambda}} F(A_\lambda') \sum_{S \in \mr{SYT}(\lambda)}  P_n \cdot e_{S'} 
= \frac{(2n)!}{f^{2\lambda}} F(A_\lambda') \sum_{T \in \mr{SYT}(2\lambda)}  P_n \cdot e_{T} 
=  F(A_\lambda') P_n \cdot \chi^{2\lambda}.
\end{align*}
Hence we obtain the desired formula because of 
$\omega^\lambda =(2^n n!)^{-1} P_n \cdot \chi^{2\lambda}$.
\end{proof}

Now we give the explicit expansion of $F(J_1,J_3,\dots,J_{2n-1}) \cdot P_n$ 
with respect to $\omega^\lambda$.

\begin{cor} \label{cor:SPexp}
For any symmetric function $F$, we have
$$
F(J_1,J_3,\dots,J_{2n-1}) \cdot P_n = 
\frac{1}{(2n-1)!!} \sum_{\lambda \vdash n} f^{2\lambda} F(A'_{\lambda}) \omega^\lambda.
$$
In particular, the multiplicity of the identity $\mr{id}_{2n}$ in $F(J_1,J_3,\dots,J_{2n-1}) \cdot P_n$
equals
\begin{equation} \label{eq:ContentEva}
\frac{1}{(2n-1)!!} \sum_{\lambda \vdash n} f^{2\lambda} F(A'_{\lambda}).
\end{equation}
\end{cor}

\begin{proof}
Recall (see (4.8) in \cite{CM})
$$
P_n =\frac{1}{(2n-1)!!} \sum_{\lambda \vdash n} f^{2\lambda} \omega^\lambda.
$$
The claim follows from this identity and Theorem \ref{thm:fOmegaContent} immediately.
\end{proof}

\section{Double coset expansion}

\subsection{Class expansion for $m_\lambda(J_1,J_2,\dots,J_n)$}
\label{subsec:ClassExp}

In this subsection,
we review some results in \cite{MN}.
These should be compared with  theorems given in the next subsection.

A permutation $\pi$ in $S_n$ is of  \emph{reduced cycle-type} $\mu$ if
$\pi$ is of the (ordinary) cycle-type $\lambda$ and 
$\mu=\tilde{\lambda}$.
Here, as defined, $\tilde{\lambda}$ is the reduction of $\lambda$.
Let $\mf{c}_\mu(n)$ be the sum of permutations in $S_{n}$ whose reduced cycle-types are $\mu$.
The element $\mf{c}_\mu(n)$ in $\mathbb{C}[S_n]$ is zero unless $|\mu|+\ell(\mu) \le n$,
and $\{\mf{c}_\mu(n) \ | \ |\mu|+\ell(\mu) \le n\}$ is a basis of the center of $\bC[S_n]$.

It is well known that,
for any $F \in \mbb{S}$,
$F(J_1,J_2,\dots,J_n)$ is an element of the center of the group algebra 
in $\bC[S_n]$,
see e.g. \cite{Jucys,MN}.
We  define coefficients $L^\lambda_\mu(n)$ 
for the monomial symmetric function $m_\lambda$
via
\begin{equation} \label{eq:defL}
m_\lambda(J_1,J_2,\dots,J_n)= \sum_{\mu:|\mu|+\ell(\mu) \le n} L^\lambda_\mu(n) \mf{c}_\mu(n).
\end{equation}

We define a number
$$
\mr{RC}(\lambda)= \frac{|\lambda|!}{(|\lambda|-\ell(\lambda)+1)! \prod_{i \ge 1} m_i(\lambda)!}.
$$
For convenience, we put $\mr{RC}(0)=1$ for the zero partition $(0)$. 
We call this the \emph{refined Catalan number}, see \cite[\S 5.1]{MN}.
It is known that $\mr{RC}(\lambda)$ is a positive integer for any $\lambda$.

\begin{thm}[\cite{MN}]   \label{thm:MN1}
Let $\lambda,\mu$ be partitions.
Then the following statements hold.
\begin{enumerate}
\item $L^\lambda_\mu(n)$ is a polynomial in $n$. 
\item $L^\lambda_\mu(n)$ is zero unless $|\lambda| \ge |\mu|$ and $|\lambda| \equiv |\mu| \pmod{2}$.
\item If $|\lambda|=|\mu|$, then $L^\lambda_\mu=L^\lambda_\mu(n)$ is independent of $n$, and given by
\begin{equation}
L^\lambda_\mu = \sum_{(\lambda^{(1)},\lambda^{(2)},\dots)} \mr{RC}(\lambda^{(1)})
\mr{RC}(\lambda^{(2)}) \cdots
\end{equation}
summed over all sequences of partitions such that
$$
\lambda^{(i)} \vdash \mu_i \ (i \ge 1) \qquad \text{and} \qquad \lambda=\lambda^{(1)} \cup \lambda^{(2)} \cup
\cdots.
$$
\end{enumerate}
\end{thm}

Define coefficients $F^k_\mu(n)$ via
\begin{equation} \label{eq:defF}
h_k(J_1,J_2,\dots,J_n)= \sum_{\mu:|\mu|+\ell(\mu) \le n} F^k_\mu(n) \mf{c}_\mu(n),
\end{equation}
where $h_k$ is the complete symmetric function of degree $k$.
Since $h_k=\sum_{\lambda \vdash k} m_\lambda$, we have
$$
F^k_\mu(n)=\sum_{\lambda \vdash k} L^\lambda_\mu(n).
$$

\begin{thm}[\cite{MN}]   \label{thm:MN2}
For $\mu \vdash k$ we have 
$$
F^k_\mu(n) = \prod_{i=1}^{\ell(\mu)} \mr{Cat}_{\mu_i}.
$$
Here $\mr{Cat}_k=\frac{1}{k+1} \binom{2k}{k}$ is the Catalan number.
\end{thm}

\subsection{Double coset expansion for $m_\lambda(J_1,J_3,\dots,J_{2n-1}) \cdot P_n$}

Like in the case of reduced cycle-types,
we prefer to reduced coset-types rather than ordinary coset-types.
A permutation $\sigma \in S_{2n}$ is of \emph{reduced coset-type} $\mu$
if $\mu=\tilde{\lambda}$ and the ordinary coset-type of $\sigma$ is $\lambda \vdash n$,
i.e., $\Xi_n(\sigma)=\lambda$.
In particular, elements in $H_n$ are of  reduced coset-type $(0)$.
If $\mu$ is the reduced coset-type of $\sigma$, we write as $\xi(\sigma)=\mu$.

Define $\psi_\mu(n)$ to be the sum of permutations in $S_{2n}$
whose reduced coset-types are $\mu$. Note that 
$\phi_\lambda= \psi_\mu(n)$ if $\lambda \vdash n$ and $\mu= \tilde{\lambda}$,
where $\phi_\lambda$ is defined in Subsection \ref{subsec:CharZonal}.
We have $\psi_\mu(n)=0$ unless $|\mu|+\ell(\mu) \le n$.
The set $\{\psi_\mu(n) \ | \ |\mu|+\ell(\mu) \le n\}$ forms a basis of the Hecke algebra $\mcal{H}_n$.

Corollary \ref{cor1-Jucys} can be rewritten as
\begin{equation}
e_k(J_1,J_3,\dots,J_{2n-1}) \cdot P_n = \sum_{\mu \vdash k} \psi_\mu(n). 
\end{equation}
We would like to generalize this formula to any monomial symmetric function $m_\lambda$.
Define coefficients $M^\lambda_\mu(n)$ by 
\begin{equation} \label{eq:defM}
m_\lambda(J_1,J_3,\dots,J_{2n-1}) \cdot P_n = \sum_\mu M^\lambda_\mu(n) \psi_\mu(n)
\end{equation}
summed over $\mu$ such that $|\mu|+\ell(\mu) \le n$.
Note that,
by Corollary \ref{cor:SPexp}, the coefficient $M^\lambda_\mu(n)$ is given by 
\begin{equation} \label{eq:Msum}
M^\lambda_\mu(n) = \frac{1}{(2n-1)!!} \sum_{\rho \vdash n} f^{2\rho} 
\omega^\rho_{\mu+(1^{n-|\mu|})} m_\lambda(A_\rho'),
\end{equation}
where $\omega^\rho_\nu$ is a value of a zonal spherical function 
defined in Subsection \ref{subsec:CharZonal}.

The following theorem is our main result for $M^\lambda_\mu(n)$.

\begin{thm} \label{thm:coefM}
Let $\lambda,\mu$ be partitions. Then
\begin{enumerate} 
\item $M^\lambda_\mu(n)$ is a polynomial in $n$. 
\item  We have the inequality
\begin{equation} \label{eq:MLinequality}
M^\lambda_\mu(n) \ge L^\lambda_\mu(n).
\end{equation}
\item $M^\lambda_\mu(n)$ is zero unless $|\lambda| \ge |\mu|$.
\item If $|\lambda|=|\mu|$, then we have $M^\lambda_\mu(n)=L^\lambda_\mu$.
Here $L^\lambda_\mu$ is given in Theorem \ref{thm:MN1}.
In particular, $M^\lambda_\mu(n)$ is independent of $n$ in this case.
\end{enumerate}
\end{thm}

Define coefficients $G^k_\mu(n)$ via
\begin{equation} \label{eq:defG}
h_k(J_1,J_3,\dots,J_{2n-1})\cdot P_n= \sum_{\mu:|\mu|+\ell(\mu) \le n} G^k_\mu(n) \psi_\mu(n),
\end{equation}
or 
$$
G^k_\mu(n)=\sum_{\lambda \vdash k} M^\lambda_\mu(n).
$$

\begin{thm} \label{thm:coefG}
For $\mu \vdash k$, we have $G^k_\mu(n)=F^k_\mu(n)=\prod_{i=1}^{\ell(\mu)} \mr{Cat}_{\mu_i}$.
\end{thm}

We will prove these theorems except 
part 1 of Theorem \ref{thm:coefM} in the coming section.
The remaining statement will be proved in Section \ref{sec:JackDeform}
by  applying the theory of shifted symmetric functions.

\section{Proof of Theorem \ref{thm:coefM} and Theorem \ref{thm:coefG}}

\subsection{Proof of part 2 of Theorem \ref{thm:coefM}}
 \label{subsec:Proof2}

Recall that quantities $L^\lambda_\mu(n)$ and $M^\lambda_\mu(n)$ are
defined by \eqref{eq:defL} and \eqref{eq:defM}, respectively.

Let $\mu$ be a  partition and let $n$ be a positive integer such that $n \ge |\mu|+\ell(\mu)$.
We define a permutation $\pi_\mu \in S_{n}$ and 
a pair partition $\mf{m}_\mu \in \mcal{M}(2n)$ as follows.
\begin{align*}
\pi_\mu=& (1 \ 2  \ \dots \ \mu_1+1)(\mu_1+2 \ \mu_1+3 \ \dots \ \mu_1+\mu_2+2) \cdots, \\
\mf{m}_\mu=&\{ \{1, 2\mu_1+2\}, \{2,3\},\dots,\{2\mu_1,2\mu_1+1\}, \\
& \quad \{2\mu_1+3,2(\mu_1+\mu_2)+4\}, \{2\mu_1+4,2\mu_1+5\},\dots,\{2(\mu_1+\mu_2)+2,2(\mu_1+\mu_2)+3\},\dots \}.
\end{align*}
For example, if $\mu=(2,1)$, we have
\begin{align*}
\pi_{(2,1)}=&(1 \ 2 \ 3)(4 \ 5)(6)(7) \cdots(n) \\
\mf{m}_{(2,1)}=&\{\{1,6\},\{2,3\},\{4,5\},\{7,10\},\{8,9\},
\{10,11\},\dots,\{2n-1,2n\} \}.
\end{align*}
By construction, the reduced cycle-type of $\sigma_\mu$ is $\mu$
and  the reduced coset-type of $\mf{m}_\mu$ is also: $\xi(\mf{m}_\mu)=\mu$.
Note that $\mf{m}_{(0)}=\{\{1,2\},\{3,4\},\dots,\{2n-1,2n\}\}$.

Define the action $\mf{L}$ of $S_{2n}$ on $\mcal{M}(2n)$ by
$$
\mf{L}(\sigma) \mf{m}= \left\{
\{\sigma(\mf{m}(1)), \sigma(\mf{m}(2))\}, \dots \dots,
\{\sigma(\mf{m}(2n-1)), \sigma(\mf{m}(2n))\} 
\right\}, \qquad (\sigma \in S_{2n}, \ \mf{m} \in \mcal{M}(2n)).
$$
Note $\mf{L}(\sigma)\mf{m}_{(0)} = \mf{m}$ if and only if $\sigma \in \mf{m} H_n$.

\begin{lem} \label{lem:NumberMdef}
Let $\lambda, \mu$ be partitions and let $n \ge |\mu|+\ell(\mu)$.
Fix $\mf{l} \in \mcal{M}(2n)$ of reduced coset-type $\mu$.
(In particular, we may take $\mf{l}=\mf{m}_\mu$.)
Then we have
$$
M^\lambda_\mu(n)= \sum_{
\begin{subarray}{c} \sigma \in S_{2n} \\ \mf{L}(\sigma) \mf{m}_{(0)}=\mf{l}
\end{subarray}} [\sigma] m_\lambda(J_1,J_3,\dots,J_{2n-1}),
$$
where $[\sigma] w$ denotes
the multiplicity of $\sigma$ in $w \in \bC[S_{2n}]$.
In particular, $M^\lambda_\mu(n)$ is a non-negative integer.
\end{lem}

\begin{proof}
From the coset decomposition \eqref{eq:rightcosets}, we have
\begin{align*}
m_\lambda(J_1,J_3,\dots,J_{2n-1}) \cdot P_n
=& \sum_{\sigma \in S_{2n}}  
([\sigma]m_\lambda(J_1,J_3,\dots,J_{2n-1})) \sigma  P_n \\
=& \sum_{ \mf{m} \in \mcal{M}(2n)} \sum_{\sigma \in \mf{m} H_n} 
([\sigma]m_\lambda(J_1,J_3,\dots,J_{2n-1})) \mf{m}  P_n \\
=& \sum_\mu \sum_{\begin{subarray}{c} \mf{m} \in \mcal{M}(2n) \\ \xi(\mf{m})=\mu
\end{subarray}}
 \sum_{\sigma \in \mf{m} H_n}  ([\sigma] m_\lambda(J_1,J_3,\dots,J_{2n-1}))\mf{m} P_n.
\end{align*}
Since
$$
\psi_\mu(n) = \sum_{\begin{subarray}{c} \sigma \in S_{2n} \\ 
\xi(\sigma)=\mu \end{subarray}}\sigma=
\sum_{\begin{subarray}{c} \mf{m} \in \mcal{M}(2n) \\ \xi(\mf{m})=\mu \end{subarray}}
\mf{m}P_n,
$$
it follows from \eqref{eq:defM} that for each $\mu$,
$$
 \sum_{\begin{subarray}{c} \mf{m} \in \mcal{M}(2n) \\ \xi(\mf{m})=\mu
\end{subarray}}
 \sum_{\sigma \in \mf{m} H_n}  ([\sigma] m_\lambda(J_1,J_3,\dots,J_{2n-1}))\mf{m} P_n
= M^\lambda_\mu(n)
\sum_{\begin{subarray}{c} \mf{m} \in \mcal{M}(2n) \\ \xi(\mf{m})=\mu \end{subarray}}
\mf{m}P_n
$$
so that
$M^\lambda_\mu(n)= \sum_{\sigma \in \mf{l} H_n}  [\sigma] m_\lambda(J_1,J_3,\dots,J_{2n-1})$. This implies the desired claim.
\end{proof}

Let $(i_1,\dots,i_k)$ be a weakly increasing sequence of $k$ positive integers.
The sequence $(i_1,\dots,i_k)$ is of \emph{type} $\lambda \vdash k$ if $\lambda=(\lambda_1,\lambda_2,\dots)$
is a permutation of $(b_1,b_2,\dots)$,
where $b_p$ is the multiplicity of $p$ in $(i_1,\dots,i_k)$.
For $\lambda \vdash k$, the monomial symmetric polynomial is expressed as 
$$
m_\lambda(x_1,x_2,\dots,x_n)= \sum_{\begin{subarray}{c} 1 \le t_1 \le t_2 \le \dots \le t_k \le n \\
(t_1,t_2, \dots,t_k): \text{type $\lambda$} \end{subarray}}
x_{t_k}\cdots x_{t_2}  x_{t_1}. 
$$

Given partitions $\lambda,\mu$ with $|\lambda|=k$, we denote by $A_n(\lambda,\mu)$
the set of 
sequences $(u_1,v_1,u_2,v_2,\dots,u_k,v_k)$ of positive integers, satisfying
the following conditions:
\begin{itemize}
\item $(v_1,v_2,\dots,v_k)$ is of type $\lambda$ and $2 \le v_1 \le v_2 \le \dots \le v_k \le n$;\item $u_i <v_i$ for all $1 \le i \le k$;
\item The product of transpositions $(u_k \ v_k) \cdots (u_2 \ v_2) (u_1 \ v_1)$ 
coincides with $\pi_\mu$.
\end{itemize}
We also denote by $B_n(\lambda,\mu)$ the set of the sequences $(s_1,t_1,\dots,s_k,t_k)$ satisfying
\begin{itemize}
\item $(t_1,t_2,\dots,t_k)$ consists of odd numbers and is
of type $\lambda$, and $3 \le t_1 \le t_2 \le \dots \le t_k \le 2n-1$;
\item $s_i <t_i$ for all $1 \le i \le k$;
\item  $\mf{L}((s_k \ t_k) \cdots (s_2 \ t_2) (s_1 \ t_1)) \mf{m}_{(0)}=\mf{m}_\mu$.
\end{itemize}
By the definitions of $L^\lambda_\mu(n)$ and Lemma \ref{lem:NumberMdef}, we have 
$$
L^\lambda_\mu(n)=|A_n(\lambda,\mu)|, \qquad M^\lambda_\mu(n)=|B_n(\lambda,\mu)|.
$$

Now the map 
$$
(u_1,v_1,u_2,v_2,\dots,u_k,v_k) \mapsto (2u_1-1, 2v_1-1, 2u_2-1,2v_2-1,\dots, 2u_k-1,2v_k-1)
$$
gives an injection from $A_n(\lambda,\mu)$ to $B_n(\lambda,\mu)$.
Indeed, suppose $(u_1,v_1,\dots,u_k,v_k)$ is an element of $A_n(\lambda,\mu)$.
Then $\sigma:=(2u_k-1 \ 2v_k-1) \cdots (2u_1-1 \ 2v_1-1)$ permutes only odd-numbered letters,
and 
$\sigma (2j-1) = 2 \pi_\mu(j)-1$ for any $j$. 
Since $\pi_\mu$ has the cycle $(1 \ 2 \ \dots \ \mu_1+1)$,
the perfect matching $\mf{L}(\sigma) \mf{m}_{(0)}$ has blocks
$\{3,2\},\{5,4\},\dots,\{2\mu_1+1,2\mu_1\}$ and $\{1,2(\mu_1+1)\}$,
which are the first $\mu_1+1$ blocks of $\mf{m}_\mu$.
Thus, we obtain $\mf{L}(\sigma) \mf{m}_{(0)}=\mf{m}_\mu$.

This injection gives $|A_n(\lambda,\mu)| \le |B_n(\lambda,\mu)|$, that is, 
$L^\lambda_\mu(n) \le M^\lambda_\mu(n)$.

\subsection{Proof of part 3 of Theorem \ref{thm:coefM}} \label{subsec:Proof3}

The discussion in this subsection is parallel to \cite[\S 5.3]{MN}.

From now, we suppose that $n$ is sufficiently large.
In Subsection \ref{subsec:PM}, we consider the inclusion $\mcal{M}(2n-2) \subset \mcal{M}(2n)$.
Under this inclusion, the reduced coset-types are invariant.

Given $\mf{m} \in \mcal{M}(2n)$, we define the set $\mcal{S}(\mf{m})$  by
$$
\mcal{S}(\mf{m}) =\big\{ j \in \{1,2,\dots,n\} | \ 
\text{$\{\mf{m}(2j-1),\mf{m}(2j)\} \not= \{2k-1,2k\}$ for any $k \ge 1$} \big\}.
$$
If the reduced coset-type of $\mf{m}$ is $\mu$, then 
$|\mcal{S}(\mf{m})|=|\mu|+\ell(\mu)$.

For a real number $x$, we put $\lceil x \rceil= \min \{n \in \bZ \ | \ x \le n \}$.
Given a positive integer $s$, define $s^{\circ}$ by
$$
s^{\circ} = \begin{cases}
s+1 & \text{if $s$ is odd}, \\
s-1 & \text{if $s$ is even}. 
\end{cases}
$$
Equivalently, $s^{\circ}$ is the unique integer satisfying 
$\{s,s^{\circ}\} \in \mf{m}_{(0)}$. We have $s=t^{\circ}$ if and only if $t=s^{\circ}$.

We use the notations in Example \ref{ex:GraphGamma}.
Given $\mf{m} \in \mcal{M}(2n)$ and  integers $1 \le i<j \le 2n$, 
the symbol $i \Lra j$ stands for $\{i,j\} \in \mf{m}$.
Also, $i \lra j$ stands for $j=i^{\circ}$.
A part of a component of the graph $\Gamma(\mf{m})$
$$
i_1 \lra i_2 \Lra i_3 \lra \dots 
$$
is called a \emph{piece} of $\Gamma(\mf{m})$.
For example, $1 \lra 2 \Lra 6 \lra 5$ is a piece of $\Gamma(\{\{1,4\},\{2,6\},\{3,5\}\})$.
If $A$ is an empty piece, the piece $ i \Lra A \Lra j $ 
stands for the piece $i \Lra  j$ simply.

\begin{lem}
Given an $\mf{m} \in \mcal{M}(2n)$ and  transposition $(s\ t)$,
let $\mf{n}=\mf{L}((s \ t)) \mf{m}$.
Suppose that $\lambda=(\lambda_1,\lambda_2,\dots)$ and 
$\mu=(\mu_1,\mu_2,\dots)$ are the reduced coset-types of $\mf{m}$ and $\mf{n}$,
respectively.
Then either $|\mu|= |\lambda|-1$, $|\mu|= |\lambda|+1$, or $\mu=\lambda$ holds.
Furthermore, if $|\mu|=|\lambda|+1$, then 
$\mcal{S}(\mf{n}) = \mcal{S}(\mf{m}) \cup \{ \lceil \tfrac{s}{2} \rceil,
\lceil \tfrac{t}{2} \rceil \}$, and
vertices $s,t$ belong to the same component of $\Gamma(\mf{n})$. 
\end{lem}

\begin{proof}
First, suppose $\lceil \tfrac{s}{2} \rceil=\lceil \tfrac{t}{2} \rceil$, i.e.,
$t=s^{\circ}$.
There exist  (possibly empty) pieces $A, B$  such that
$A \Lra s \lra t \Lra B$ is a piece of $\Gamma(\mf{m})$,
and then $\Gamma(\mf{n})$ has the piece 
$A \Leftrightarrow t \leftrightarrow s \Leftrightarrow B$.
Therefore we have $\lambda=\mu$.

From now, we suppose $\lceil \tfrac{s}{2} \rceil\not=\lceil \tfrac{t}{2} \rceil$,
and so $s,s^{\circ},t$, and $t^{\circ}$ are distinct.
Then the following five cases may occur:
\begin{itemize}
\item[(i)] $|\mcal{S}(\mf{m}) \cap \{\lceil \tfrac{s}{2} \rceil, \lceil \tfrac{t}{2} \rceil \}|=0$.
\item[(ii)] $|\mcal{S}(\mf{m}) \cap \{\lceil \tfrac{s}{2} \rceil, \lceil \tfrac{t}{2} \rceil \}|=1$.
\item[(iii)] $\{\lceil \tfrac{s}{2} \rceil, \lceil \tfrac{t}{2} \rceil \} \subset \mcal{S}(\mf{m})$
and $s,t$ belong to different components of $\Gamma(\mf{m})$.
\item[(iv)] $\Gamma(\mf{m})$ has a component of the form
\begin{equation} \label{eq:componentst1}
s \lra s^{\circ} \Lra A \Lra t \lra t^{\circ} \Lra  B \Lra s. 
\end{equation}
\item[(v)] $\Gamma(\mf{m})$ has a component of the form
\begin{equation} \label{eq:componentst2}
s \lra s^{\circ} \Lra C \Lra t^{\circ} \lra t \Lra  D \Lra s.
\end{equation}
\end{itemize}
Here $A, B, C$ and $D$ are possibly empty pieces.
For each case, we shall compare $\Gamma(\mf{n})=\Gamma(\mcal{L}(s \ t) \mf{m})$ 
with $\Gamma(\mf{m})$.

In the case (i), 
the graph $\Gamma(\mf{m})$ has components $s\lra s^{\circ} \Lra s$ and
$t \lra t^{\circ} \Lra t$, and
the graph $\Gamma(\mf{n})$ has the new component
$s \lra s^{\circ} \Lra t \lra t^{\circ} \Lra s$. 
Thus $\mu = \lambda \cup (1)$.

In the case (ii), we may suppose $\mcal{S}(\mf{m}) \cap \{\lceil \tfrac{s}{2} \rceil, \lceil \tfrac{t}{2} \rceil \}=\{ \lceil \tfrac{s}{2} \rceil\}$.
A piece $A \Lra s \lra s^{\circ} \Lra B$ of $\Gamma(\mf{m})$ with some pieces $A,B$
becomes
the piece  $A \Lra t \lra t^{\circ} \Lra s  \lra s^{\circ} \Lra B$ of $\Gamma(\mf{n})$.
Therefore $\mu$ has a part equal to $\lambda_j+1$. 
In particular, $\mcal{S}(\mf{n}) = \mcal{S}(\mf{m}) \cup \{ \lceil \tfrac{s}{2} \rceil,
\lceil \tfrac{t}{2} \rceil \}$.

In case (iii), $\Gamma(\mf{m})$ has two components of the forms
$s \lra s^{\circ} \Lra A \Lra s$ and
$t \lra t^{\circ} \Lra B \Lra t$,
where $A,B$ are non-empty pieces.
Then $\Gamma(\mf{n})$ has the combined component 
$$
s \lra s^{\circ} \Lra A \Lra t \lra t^{\circ} \Lra B \Lra s.
$$
Therefore a certain part $\mu_m$ of $\mu$ equals $\lambda_i+\lambda_j+1$ 
for some $1 \le i<j \le \ell(\lambda)$. 
We also see that $\{\lceil \tfrac{s}{2} \rceil, \lceil \tfrac{t}{2} \rceil \} \subset \mcal{S}(\mf{m}) =\mcal{S}(\mf{n})$.

In case (iv),  $\Gamma(\mf{n})$ has divided components 
$$
s \lra s^{\circ} \Lra A \Lra s \qquad \text{and} \qquad
t \lra t^{\circ} \Lra B \Lra t.
$$ 
Therefore there are $\mu_i$ and $\mu_j$ equal to $r-1$ and $\lambda_m-r$
for some $\lambda_m$ and $1 \le r \le \lambda_m$.
In particular, $|\mu|=|\lambda|-1$.

In case (v),
$\Gamma(\mf{n})$ has the component 
$$
s \lra s^{\circ} \Lra C \Lra t^{\circ} \lra t \Lra D^\vee \Lra s.
$$ 
Here, if $D$ is the piece $i_1 \lra i_2 \Lra \cdots \lra i_{2p}$ then
 $D^\vee$ is the piece $i_{2p} \lra \cdots \Lra i_2 \lra i_1$.
In this case, $\lambda=\mu$.

For the only cases (i), (ii), and (iii), we have $|\mu|=|\lambda|+1$.
The rest of the claims are already seen.
\end{proof}

\begin{cor} \label{cor:numbertrans}
Let $\mu$ be the reduced coset-type of $\mf{m} \in \mcal{M}(2n)$.
Suppose that there exist $p$ transpositions $(s_1 \ t_1), \dots, (s_p \ t_p)$
satisfying
$\mf{L}( (s_p\ t_p) \cdots (s_1 \ t_1) ) \mf{m}_{(0)} = \mf{m}$.  
Then $|\mu| \le p$.
\end{cor}

\begin{cor} \label{cor:BestTrans}
Let $\mu \vdash p$ and let $\mf{m} \in \mcal{M}(2n)$ be of reduced coset-type $\mu$.
Suppose that there exist $p$ transpositions $(s_1 \ t_1), \dots, (s_p \ t_p)$
satisfying
$\mf{L}( (s_p \ t_p) \cdots (s_1 \ t_1) ) \mf{m}_{(0)} = \mf{m}$.  
Then $\mcal{S}(\mf{m})=\{
\lceil \tfrac{s_1}{2} \rceil, \lceil \tfrac{t_1}{2} \rceil,
\dots, \lceil \tfrac{s_p}{2} \rceil, \lceil \tfrac{t_p}{2} \rceil \}$.
Furthermore, for each $i$, the vertices $s_i,t_i$ belong to the same component of 
$\Gamma(\mf{m})$.
\end{cor}

Since $m_\lambda(J_1,J_3,\dots,J_{2n-1})$ is 
a sum of products of $|\lambda|$ transpositions,
part 3 of Theorem \ref{thm:coefM} follows from Corollary \ref{cor:numbertrans}
together with Lemma \ref{lem:NumberMdef}.

\subsection{Proof of Theorem \ref{thm:coefG}}

Recall that quantities $F^k_\mu(n)$ and $G^k_\mu(n)$ are
defined by \eqref{eq:defF} and \eqref{eq:defG}, respectively.

Let $\mf{m}, \mf{n} \in \mcal{M}(2n)$ and 
suppose $\mcal{S}(\mf{m}) \cap \mcal{S}(\mf{n}) =\emptyset$.
Denote by $\tilde{\mf{m}}$  the perfect matching on 
$\bigsqcup_{i \in \mcal{S}(\mf{m})} \{2i-1,2i\}$ obtained as
the union of 
non-trivial blocks of $\mf{m}$.
Clearly, $\mcal{S}(\mf{m})=\mcal{S}(\tilde{\mf{m}})$.
We define the new perfect matching $\mf{m} \cup \mf{n} \in \mcal{M}(2n)$  by
$$
\mf{m} \cup \mf{n} = \tilde{\mf{m}} \sqcup \tilde{\mf{n}} \sqcup
\{\{ 2i-1,2i \} \ | \ i \not\in \mcal{S}(\mf{m}) \sqcup \mcal{S}(\mf{n}) \}.
$$
If $\lambda$ and $\mu$ is the reduced coset-type of $\mf{m}$ and $\mf{n}$, respectively,
then the reduced coset-type of $\mf{m} \cup \mf{n}$ is $\lambda \cup \mu$.

\begin{example}
For
\begin{align*}
\mf{m}=&\{ \{1,5\},\{3,4\},\{2,6\},  \{7,8\},\{9,10\},\{11,12\},\dots,\{2n-1,2n\}\}, \\
\mf{n}=&\{\{1,2\},\{3,4\},\{5,6\}, \{7,10\},\{8,9\},\{11,12\},\dots,\{2n-1,2n\}\},
\end{align*}
we have
$$
\mf{m} \cup \mf{n}= \{\{1,5\},\{2,6\},\{3,4\},\{7,10\},\{8,9\},\{11,12\},\dots,\{2n-1,2n\} \}.
$$
The reduced coset-types of $\mf{m}$, $\mf{n}$, and $\mf{m} \cup \mf{n}$  
are $(1)$, $(1)$, and $(1,1)$, respectively.
\end{example}

\begin{lem} \label{lem:PMdecomCoe}
Let $\mf{n}^{(1)}, \mf{n}^{(2)}\in \mcal{M}(2n)$
such that $k<l$ for all $k \in \mcal{S}(\mf{n}^{(1)})$ and
$l \in \mcal{S}(\mf{n}^{(2)})$.
Suppose that the reduced coset-types of $\mf{n}^{(i)}$ have 
sizes $r_i$ \ ($i=1,2$).
Also, there exist $r$ transpositions $(s_1 \ t_1),\dots, (s_r \ t_r)$
satisfying
$\mf{L}((s_r \ t_r) \cdots (s_1 \  t_1)) \mf{m}_{(0)} =\mf{n}^{(1)} \cup \mf{n}^{(2)}$,
where $r=r_1+r_2, s_i<t_i \ (1 \le i \le r)$, and 
$t_r \ge \cdots \ge t_1$.
Then 
$$
\mf{n}^{(1)}= \mf{L}( (s_{r_1} \ t_{r_1}) \cdots (s_1 \ t_1)) \mf{m}_{(0)},
\qquad \mf{n}^{(2)}= \mf{L}( (s_{r} \ t_{r}) \cdots (s_{r_1+1} \ t_{r_1+1})) \mf{m}_{(0)}
$$
and 
$$
\mcal{S}(\mf{n}^{(1)}) = \{
\lceil \tfrac{s_1}{2} \rceil, \lceil \tfrac{t_1}{2} \rceil,
\dots, \lceil \tfrac{s_{r_1}}{2} \rceil, \lceil \tfrac{t_{r_1}}{2} \rceil \},
\qquad
\mcal{S}(\mf{n}^{(2)})=\{
\lceil \tfrac{s_{r_1+1}}{2} \rceil, \lceil \tfrac{t_{r_1+1}}{2} \rceil,
\dots, \lceil \tfrac{s_r}{2} \rceil, \lceil \tfrac{t_r}{2} \rceil \}.
$$
\end{lem}

\begin{proof}
By Corollary \ref{cor:BestTrans}, we see
$\mcal{S}(\mf{n}^{(1)}) \sqcup \mcal{S}(\mf{n}^{(2)}) =
\mcal{S}(\mf{n}^{(1)} \cup \mf{n}^{(2)})=\{
\lceil \tfrac{s_1}{2} \rceil, \lceil \tfrac{t_1}{2} \rceil,
\dots, \lceil \tfrac{s_r}{2} \rceil, \lceil \tfrac{t_r}{2} \rceil \}$.
Since $t_i$ are not decreasing, there exists an integer $p$ such that
$\{\lceil \tfrac{t_1}{2} \rceil,
\dots, \lceil \tfrac{t_p}{2} \rceil \} \subset \mcal{S}(\mf{n}^{(1)})$
and $\{\lceil \tfrac{t_{p+1}}{2} \rceil,
\dots, \lceil \tfrac{t_r}{2} \rceil \} \subset \mcal{S}(\mf{n}^{(2)})$.
Furthermore, applying Corollary \ref{cor:BestTrans} again,
we see that $s_i,t_i$ belong to the same component of 
$\Gamma(\mf{n}^{(1)} \cup \mf{n}^{(2)})$,
and so that $\mcal{S}(\mf{n}^{(1)})=\{
\lceil \tfrac{s_1}{2} \rceil, \lceil \tfrac{t_1}{2} \rceil,
\dots, \lceil \tfrac{s_p}{2} \rceil, \lceil \tfrac{t_p}{2} \rceil \}$
and 
$\mcal{S}(\mf{n}^{(2)})=\{
\lceil \tfrac{s_{p+1}}{2} \rceil, \lceil \tfrac{t_{p+1}}{2} \rceil,
\dots, \lceil \tfrac{s_r}{2} \rceil, \lceil \tfrac{t_r}{2} \rceil \}$.

Put $\rho^{(1)}= (s_p \ t_p) \cdots (s_1 \ t_1)$ and 
$\rho^{(2)}=  (s_r \ t_r) \cdots (s_{p+1} \ t_{p+1})$.
Since $\{s_1,t_1,\dots,s_p,t_p\} \cap \{s_{p+1},t_{p+1},\dots,s_r,t_r\}=\emptyset$,
we have
 $\mf{n}^{(1)} \cup \mf{n}^{(2)} =\mcal{L}(\rho^{(2)} \rho^{(1)})\mf{m}_{(0)} =
\mcal{L}(\rho^{(1)})\mf{m}_{(0)}
\cup \mcal{L}(\rho^{(2)})\mf{m}_{(0)}$
and
\begin{align*}
\mcal{S}(\mf{L}(\rho^{(1)})\mf{m}_{(0)}) =& \{
\lceil \tfrac{s_1}{2} \rceil, \lceil \tfrac{t_1}{2} \rceil,
\dots, \lceil \tfrac{s_p}{2} \rceil, \lceil \tfrac{t_p}{2} \rceil \}
=\mcal{S}(\mf{n}^{(1)}), \\
\mcal{S}(\mf{L}(\rho^{(2)})\mf{m}_{(0)}) =& \{
\lceil \tfrac{s_{p+1}}{2} \rceil, \lceil \tfrac{t_{p+1}}{2} \rceil,
\dots, \lceil \tfrac{s_r}{2} \rceil, \lceil \tfrac{t_r}{2} \rceil \}
=\mcal{S}(\mf{n}^{(2)}).
\end{align*}
Hence $\mf{n}^{(i)}= \mf{L}(\rho^{(i)})\mf{m}_{(0)}$ \ $(i=1,2)$.
In particular, the size of the reduced coset-type of $\mf{L}(\rho^{(i)})\mf{m}_{(0)}$
is $r_i$ \ $(i=1,2)$.
On the other hand, Corollary \ref{cor:numbertrans} and the definition of $\rho^{(i)}$ imply that 
$r_1 \le p$ and $r_2 \le r-p$.
But $r=r_1+r_2$ so that $p=r_1$.
\end{proof}

Given a positive integer $k$ and a perfect matching $\mf{l} \in \mcal{M}(2n)$,
we define $\mcal{B}_n(k,\mf{l})$ by 
the set of all sequences $(s_1,t_1,s_2,t_2,\dots,s_k,t_k)$ of positive integers,
satisfying
following conditions.
\begin{itemize}
\item All of $t_i$ are odd and $3 \le t_1 \le \cdots \le t_{k} \le 2n-1$;
\item $s_i <t_i$ for all $i$;
\item $\mf{L}((s_k \ t_k) \cdots (s_1 \ t_1))\mf{m}_{(0)}= \mf{l}$.
\end{itemize}
Remark that the set $\mcal{B}_n(k,\mf{m}_\mu)$ coincides with the union 
$\bigsqcup_{\lambda \vdash k} B_n(\lambda,\mu)$,
where $B_n(\lambda,\mu)$ was defined in Subsection \ref{subsec:Proof2}.
Also we define the set $\mcal{B}(k,\mf{l})$ 
as the subset of $\mcal{B}_n(k,\mf{l})$ which consists of
sequences satisfying
$$
\{\lceil \tfrac{t_1}{2} \rceil, \dots, \lceil \tfrac{t_k}{2} \rceil \} \subset \mcal{S}(\mf{l})
\qquad \text{and} \qquad t_k= 2 a-1,
$$
where $a$ is the maximum of $\mcal{S}(\mf{l})$.

\begin{lem} \label{lem:coeGn}
Let $\mu \vdash k$ and let $\mf{l}$ be a perfect matching of reduced coset-type $\mu$.
Then $\mcal{B}_n(k,\mf{l}) = \mcal{B}(k,\mf{l})$ and 
$G^k_\mu(n)= |\mcal{B}(k,\mf{l})|$.
In particular, both $\mcal{B}_n(|\mu|,\mf{l})$ and $G^k_\mu(n)$ are independent of $n$.
\end{lem}

\begin{proof}
Let $(s_1,t_1,\dots,s_k,t_k)$ be an element in $\mcal{B}_n(k,\mf{l})$.
Then by Corollary \ref{cor:BestTrans}, we have
$$
\{\lceil \tfrac{s_1}{2} \rceil, \lceil \tfrac{t_1}{2} \rceil,
\dots, \lceil \tfrac{s_k}{2} \rceil, \lceil \tfrac{t_k}{2} \rceil 
\}=\mcal{S}(\mf{l}).
$$
Therefore $t_k=2 a-1$ with $a= \max \mcal{S}(\mf{l})$. 
Hence $(s_1,t_1,\dots,s_k,t_k) \in \mcal{B}(k,\mf{l})$, and so
$\mcal{B}(k,\mf{l})=\mcal{B}_n(k,\mf{l})$.
Also we have $G^k_\mu(n)=|\mcal{B}_n(k,\mf{l})|$
from Lemma \ref{lem:NumberMdef}.
\end{proof}

\begin{lem}  \label{lem:Gdecomp}
Let $\mu \vdash k$.
Then 
$G^k_\mu(n)= \prod_{i=1}^{\ell(\mu)} G^{\mu_i}_{(\mu_i)}(n)$.
\end{lem}

\begin{proof}
We prove by induction on $\ell(\mu)=l$.
If $l=1$, then our claim is trivial. Assume $l>1$.

The perfect matching $\mf{m}_\mu$ may be uniquely expressed as $\mf{m}_\mu=\mf{m}_{\nu}
\cup \mf{n}$,
where $\nu=(\mu_1,\mu_2,\dots,\mu_{l-1})$, and 
$\mf{n}$ is the perfect matching such that 
$\mcal{S}(\mf{n})=\{\mu_1+\cdots+\mu_{l-1}+l+j \ | \ 0 \le j \le \mu_l\}$.

Let $(s_1,t_1,\dots,s_k,t_k)$ be a sequence in $\mcal{B}(k,\mf{m}_\mu)$.
By Lemma \ref{lem:PMdecomCoe},
this sequence  satisfies 
$$
\mf{L}( (s_{k-\mu_l} \ t_{k-\mu_l}) \dots (s_1 \ t_1)) = \mf{m}_{\nu}, \qquad
\mf{L}( (s_{k} \ t_{k}) \dots (s_{k-\mu_l+1} \ t_{k-\mu_l+1})) = \mf{n},
$$
and 
\begin{gather*}
3 \le  t_1 \le \cdots \le t_{k-\mu_l} = 2(\mu_1+\cdots+\mu_{l-1}+l-1)-1, \\
 2(\mu_1+\cdots+\mu_{l-1}+l)+1 \le  t_{k-\mu_l+1} \le \cdots \le t_k = 2(k+l)-1.
\end{gather*}
Therefore $(s_1,t_1,\dots,s_{k-\mu_l},t_{k-\mu_l})$ belongs to $\mcal{B}(k-\mu_l,\mf{m}_\nu)$, and
$(s_{k-\mu_l+1}, t_{k-\mu_l+1}, \dots, t_{k})$
belongs to $\mcal{B}(\mu_l,\mf{n})$.
This gives a bijection between 
$\mcal{B}(k,\mf{m}_\mu)$ and $\mcal{B}(k-\mu_l,\mf{m}_\nu) \times \mcal{B} (\mu_l,\mf{n})$.
Hence the claim follows from Lemma \ref{lem:coeGn} and the assumption of the induction.
\end{proof}

\begin{lem}
$G^{k}_{(k)}(n)= \mr{Cat}_k$.
\end{lem}

\begin{proof}
We prove by induction on $k$.
Assume that for any $0 \le q <k$ it holds that $G^q_{(q)}(n)=\mr{Cat}_q$.
Let $(s_1,t_1,\dots,s_k,t_k)$ be an element of $\mcal{B}(k,\mf{m}_{(k)})$.
Then $t_k=2k+1$.
Put $p=s_k$ and $\mf{n}=\mf{L} ((s_{k-1} \ t_{k-1}) \cdots (s_1 \ t_1)) \mf{m}_{(0)}$.
Note that $\mf{n}= \mf{L}( (p \ 2k+1))\mf{m}_{(k)}$
and that the reduced coset-type of $\mf{n}$ 
is of size $k-1$.

Suppose $p$ is even, say $p=2q$ with $1 \le q \le k$.
Since $\mf{m}_{(k)}=\{\{1,2k+2\},\{2,3\},\dots,\{2k,2k+1\}\}$,
the graph $\Gamma( \mf{L}( (2q \ 2k+1))\mf{m}_{(k)})$ has 
only one non-trivial component 
$$
1 \lra 2 \Lra 3 \lra  4 \Lra \dots \lra 2q \Lra 2k \lra 2k-1 \Lra 2k-2 \lra \dots
\lra 2q+1 \Lra 2k+1 \lra 2k+2 \Lra 1.
$$
Therefore the reduced coset-type of $\mf{n}$ is $(k)$ but this is contradictory.
Hence $p=s_k$ must be odd.

Write as $s_k=p=2q-1$ with $1 \le q \le k$.
The perfect matching $\mf{n}$ can be expressed as $\mf{n}=\mf{n}_q' \cup \mf{n}_q''$, 
where $\mf{n}_q'$ and $\mf{n}_q''$ are perfect matchings in $\mcal{M}(2n)$ such that
\begin{align*}
\widetilde{\mf{n}_q'}=& \{\{1,2k+2\},\{2,3\},\{4,5\},\dots,\{2q-4,2q-3\},\{2q-2,2k+1\}\}, \\
\widetilde{\mf{n}_q''}=& \{\{2q-1,2k\},\{2q,2q+1\}, \{2q+2,2q+3\}, \dots,\{2k-2,2k-1\}  \}.
\end{align*}
The reduced coset-type of $\mf{n}$ is either $(q-1,k-q)$ or $(k-q,q-1)$.
Therefore the sequence $(s_1,t_1,\dots, s_{k-1},t_{k-1})$ belongs to $\mcal{B}_n(k-1,\mf{n})$.
Conversely, if $(s_1',t_1',\dots,s_{k-1}',t_{k-1}')$ is an element of 
$\mcal{B}_n(k-1,\mf{n}_q' \cup \mf{n}_q'')$,
then $(s_1',t_1',\dots,s_{k-1}',t_{k-1}', 2q-1, 2k+1)$ belongs to $\mcal{B}_n(k,\mf{m})$.
Therefore we have the identity
$$
G_{(k)}^k(n)=|\mcal{B}_n(k,\mf{m})|= \sum_{q=1}^k |\mcal{B}_n(k-1,\mf{n}_q' \cup \mf{n}_q'')|.
$$
It follows from Lemma \ref{lem:coeGn}, Lemma \ref{lem:Gdecomp}, and the induction assumption that
$$
 |\mcal{B}_n(k-1,\mf{n}_q' \cup \mf{n}_q'')| =
 G^{q-1}_{(q-1)}(n)G^{k-q}_{(k-q)}(n) = 
 \mr{Cat}_{q-1} \mr{Cat}_{k-q}.
$$
Hence the well known recurrence formula for Catalan numbers gives
$G_{(k)}^k(n)=\sum_{q=1}^k \mr{Cat}_{q-1} \mr{Cat}_{k-q}=\mr{Cat}_k$.
\end{proof}

We have obtained the proof of Theorem \ref{thm:coefG}.

\subsection{Proof of part 4 of Theorem \ref{thm:coefM}} 

Let $\mu \vdash k$.
By Theorem \ref{thm:MN2}, part 2 of Theorem \ref{thm:coefM}, and Theorem \ref{thm:coefG},
we see
$$
\prod_{i=1}^{\ell(\mu)}\mr{Cat}_{\mu_i}
=F^k_\mu(n)= \sum_{\lambda \vdash k} L^\lambda_\mu(n) \le 
\sum_{\lambda \vdash k} M^\lambda_\mu(n) = G^k_\mu(n) = \prod_{i=1}^{\ell(\mu)}\mr{Cat}_{\mu_i}
$$
so that $M^\lambda_\mu(n)= L^\lambda_\mu(n)$ for all $\lambda \vdash k$.

\section{Weingarten functions for the orthogonal group}
\label{section:WgOrtho}

Fix  positive integers $N,n$ and assume $N \ge n$.
We define the Weingarten function for the orthogonal group $O(N)$ by
\begin{equation} \label{eq:WgDefinition}
\mr{Wg}^{O(N)}_{n}= \frac{1}{(2n-1)!!} \sum_{\lambda \vdash n} \frac{f^{2\lambda}}
{\prod_{\square \in \lambda} (N+c'(\square))} \omega^\lambda,
\end{equation}
which is an element of the Hecke algebra $\mcal{H}_n$ of the Gelfand pair $(S_{2n},H_n)$.
Here $f^{2\lambda}$ and $c'(\square)$ were defined in Section \ref{Sec:Preparations}.
As proved in \cite{CM},
this Weingarten function plays an important role in  calculations of integrals
of polynomial functions over the orthogonal group $O(N)$.

\begin{prop}[\cite{CM}, see also \cite{CS}] \label{prop:CM}
Suppose $N \ge n$.
Let $g=(g_{ij})_{1 \le i,j \le N}$ be a Haar-distributed random matrix from $O(N)$ and
let $d g$ the normalized Haar measure on $O(N)$. 
Given two functions $\bm{i},\bm{j}$ from $\{1,2,\dots,2n\}$ to $\{1,2,\dots, N\}$,
we have 
\begin{align*}
& \int_{g \in O(N)} g_{\bm{i}(1) \bm{j}(1)} g_{\bm{i}(2) \bm{j}(2)} \cdots  g_{\bm{i}(2n) \bm{j}(2n)} d g \\
=& \sum_{\mf{m},\mf{n} \in \mcal{M}(2n)} \mr{Wg}^{O(N)}_{n}(\mf{m}^{-1} \mf{n})
\prod_{k=1}^n \delta_{\bm{i}(\mf{m}(2k-1)),\bm{i}(\mf{m}(2k))} 
\delta_{\bm{j}(\mf{m}(2k-1)),\bm{j}(\mf{m}(2k))}.
\end{align*}
Here we regard $\mcal{M}(2n)$ as a subset of $S_{2n}$.
\end{prop}

As a special case of Proposition \ref{prop:CM}, we obtain an integral expression for $\mr{Wg}^{O(N)}_n(\sigma)$:
$$
\mr{Wg}^{O(N)}_n(\sigma) = \int_{g \in O(N)} g_{1 j_1} g_{1 j_2} g_{2 j_3} g_{2 j_4} 
\cdots g_{n j_{2n-1}} g_{n j_{2n}} dg, \qquad 
\sigma \in S_{2n},
$$
with
$$
(j_1,j_2,\dots,j_{2n})= \(\left\lceil \tfrac{\sigma(1)}{2} \right\rceil, 
\left\lceil \tfrac{\sigma(2)}{2} \right\rceil, \dots, \left\lceil \tfrac{\sigma(2n)}{2} \right\rceil \).
$$

\begin{remark}
In Proposition \ref{prop:CM}, we can remove the assumption $N \ge n$.
In fact, when $N <n$, it is enough to replace the range of the sum on \eqref{eq:WgDefinition}
by partitions $\lambda \vdash n$ such that $\ell(\lambda) \le N $. 
See \cite{CM} for details.
\end{remark}

Recall that the generating function for complete symmetric polynomials $h_k$ is 
$$
\sum_{k=0}^\infty h_k(x_1,x_2,\dots,x_n) u^k = \prod_{i =1}^n \frac{1}{1-x_i u}.
$$

\begin{thm} \label{thm:OrthoWgJM}
Suppose $N \ge 2n-1$. Then
$$
\mr{Wg}^{O(N)}_{n}= \sum_{k =0}^\infty (-1)^k N^{-n-k}  h_k(J_1,J_3,\dots,J_{2n-1}) \cdot P_n.
$$
\end{thm}

\begin{proof}
We have
\begin{align*}
\mr{Wg}_n^{O(N)} =& \frac{1}{(2n-1)!!}\sum_{\lambda \vdash n}
f^{2\lambda} \(\prod_{\square \in \lambda} (N+c'(\square))^{-1} \) \omega^\lambda \\
=& \frac{1}{(2n-1)!!}\sum_{\lambda \vdash n}
f^{2\lambda} \(\sum_{k =0}^\infty (-1)^k N^{-n-k}  h_k(A_\lambda') \) \omega^\lambda \\
=&  \frac{1}{(2n-1)!!}\sum_{k =0}^\infty (-1)^k N^{-n-k} \sum_{\lambda \vdash n}
f^{2\lambda} h_k(A_\lambda') \omega^\lambda \\
=& \sum_{k =0}^\infty(-1)^k N^{-n-k}  h_k(J_1,J_3,\dots,J_{2n-1}) \cdot P_n. 
\end{align*}
Here the second equality follows because of $|c'(\square)| \le 2n-2 <N$
for all $\square \in \lambda \vdash n$, and the fourth equality follows by Corollary \ref{cor:SPexp}. 
\end{proof}

Recall that  $G_\mu^k(n)$ are coefficients in
$$
h_k(J_1,J_3,\dots,J_{2n-1}) \cdot P_n = \sum_\mu G_\mu^k (n) \psi_\mu(n).
$$
These coefficients appear in the asymptotic expansion of
$\mr{Wg}^{O(N)}_n(\sigma)$ with respect to $\tfrac{1}{N}$. 

\begin{thm} \label{thm:ExpansionWg}
Let $\mu$ be a partition and let $N,n,k$ be positive integers. Suppose 
$N \ge 2n-1$ and $n \ge |\mu|+\ell(\mu)$.
For any permutation $\sigma$ in $S_{2n}$ of reduced coset-type $\mu$,
we have
\begin{align} 
\mr{Wg}^{O(N)}_{n}(\sigma)=& 
\sum_{g=0}^\infty (-1)^{|\mu|+g} G_\mu^{|\mu|+g}(n) N^{-n-|\mu|-g} \label{eq:AsymExpWg}\\
=& (-1)^{|\mu|} \prod_{i = 1}^{\ell(\mu)} \mr{Cat}_{\mu_i} \cdot N^{-n-|\mu|}
+ (-1)^{|\mu|+1} \prod_{i = 1}^{\ell(\mu)} G^{|\mu|+1}_\mu(n) \cdot N^{-n-|\mu|-1}
+\cdots.
\end{align}
\end{thm}

\begin{proof}
Theorem \ref{thm:OrthoWgJM} and the definition of $G^k_\mu(n)$ imply
$$
\mr{Wg}^{O(N)}_{n}(\sigma)= \sum_{k=0}^\infty (-1)^k G_\mu^k(n) N^{-n-k}.
$$
It follows from Theorem \ref{thm:coefM} and Theorem \ref{thm:coefG} that
$G^{k}_\mu(n)=\sum_{\lambda \vdash k} M^\lambda_\mu(n)$ is zero unless $k \ge |\mu|$ and
that $G^{|\mu|}_{\mu}(n)=\prod_{i=1}^{\ell(\mu)} \mr{Cat}_{\mu_i}$. 
\end{proof}

The unitary group version of results in this section is seen in \cite{MN}.

Collins and \'{S}niady \cite{CS} obtained 
$$
\mr{Wg}^{O(N)}_{n}(\sigma)= (-1)^{|\mu|} \prod_{i \ge 1} \mr{Cat}_{\mu_i} \cdot N^{-n-|\mu|} +
\mr{O} (N^{-n-|\mu|-1}), \qquad N \to \infty,
$$
where $\sigma$ is a permutation in $S_{2n}$ of reduced coset-type $\mu$.
Our result is a refinement of their one.

We will observe the subleading coefficient $G^{|\mu|+1}_\mu(n)$ later,
see Subsection \ref{subsec:Open}.

\section{Jack deformations}  \label{sec:JackDeform}

A purpose in this section is to intertwine $M^\lambda_\mu(n)$ with
$L^\lambda_\mu(n)$.
They have been defined via symmetric functions in Jucys-Murphy elements.
We define their $\alpha$-extension based on the theory of Jack polynomials.

\subsection{Jack-Plancherel measures}

Let $\alpha>0$ be a positive real number.

For each $\lambda \vdash n$, we put
$$
j_{\lambda}^{(\alpha)}= \prod_{(i,j) \in \lambda} \left\{(\alpha(\lambda_i-j) + \lambda_j'-i +1)
(\alpha(\lambda_i-j) + \lambda_j'-i +\alpha) \right\},
$$
where $\lambda'=(\lambda_1',\lambda_2',\dots)$ is the conjugate partition of $\lambda$.
Here the Young diagram of $\lambda'$ is, by definition, 
the transpose of the Young diagram $\lambda$.
Define
\begin{equation}
\PP_n^{(\alpha)} (\lambda)= \frac{\alpha^n n!}{j_\lambda^{(\alpha)}}.
\end{equation}
This gives a probability measure on partitions of $n$ and is called the {\it Jack-Plancherel measure} or \emph{Jack measure} shortly.
When $\alpha=1$, 
\begin{equation}
\PP_n^{(1)} (\lambda)= \frac{n!}{(H_\lambda)^2} = \frac{(f^\lambda)^2}{n!},
\end{equation}
where $H_\lambda= \sqrt{j_\lambda^{(1)}}$ is the product of hook-lengths of $\lambda$,
and the well-known hook-length formula gives $f^\lambda = \frac{n!}{H_\lambda}$.
The probability measure $\PP^{(1)}_n$ is known as the {\it Plancherel measure}
for the symmetric group $S_n$.
Also, it is easy to see that
$$
\PP_n^{(2)} (\lambda)= \frac{f^{2\lambda}}{(2n-1)!!}, \qquad
\PP^{(1/2)}_n(\lambda)= \frac{f^{\lambda \cup \lambda}}{(2n-1)!!}. 
$$

\begin{example}
$$
\PP_3^{(\alpha)} ((3))=\frac{1}{(1+\alpha)(1+2\alpha)}, \quad
\PP_3^{(\alpha)} ((2,1))=\frac{6\alpha}{(2+\alpha)(1+2\alpha)}, \quad 
\PP_3^{(\alpha)} ((1^3))=\frac{\alpha^2}{(1+\alpha)(2+\alpha)}.
$$
\end{example}

The Jack-Plancherel measure has the duality relation:
$$
\PP_n^{(\alpha)}(\lambda) = \PP_n^{(\alpha^{-1})}(\lambda'),
$$
which follows from $j_\lambda^{(\alpha)}= \alpha^{2|\lambda|} j_{\lambda'}^{(\alpha^{-1})}$.

Some asymptotic properties of random variables $\lambda_1,\lambda_2,\dots$
with repect to Jack-Plancherel measures in $n \to \infty$ are studied, 
see  \cite{Mat} and its references.

\subsection{Jack symmetric functions}

Recall the fundamental properties for Jack symmetric functions $J^{(\alpha)}_\lambda$.
The details are seen in \cite[VI-10]{Mac}.

Consider a scalar product on 
the algebra $\mathbb{S}$ of symmetric functions given by
$$
\langle p_\lambda, p_\mu \rangle_\alpha = \delta_{\lambda, \mu} \alpha^{\ell(\lambda)}
z_\lambda
$$
where $z_\lambda$ is defined in \eqref{eq:zlambda}.
The Jack functions $\{J_\lambda^{(\alpha)} \ | \ \text{$\lambda$: partitions} \}$ are
the unique family satisfying the following two conditions:
\begin{itemize}
\item They are of the form
$J_\lambda^{(\alpha)} = \sum_{\mu \le \lambda} u_{\lambda \mu}^{(\alpha)} m_\mu$,
where each coefficient $u^{(\alpha)}_{\lambda \mu}$ is a rational function in $\alpha$,
 and where $\mu \le \lambda$ stands for
the dominance ordering: $|\mu|=|\lambda|$ and $\mu_1+\cdots+\mu_i \le \lambda_1+\cdots+\lambda_i$
for any $i \ge 1$.  
\item (orthogonality) $\langle J_\lambda^{(\alpha)}, J^{(\alpha)}_\mu \rangle_\alpha= \delta_{\lambda,\mu}
j_\lambda^{(\alpha)}$ for any $\lambda,\mu$.
\end{itemize}
We note $J^{(1)}_\lambda=H_\lambda s_\lambda$ and $J^{(2)}_\lambda=Z_\lambda$,
where $s_\lambda$ is a Schur function and $Z_\lambda$ is a zonal polynomial.

Let $\theta^\lambda_\rho(\alpha)$ be the coefficient of $p_\rho$ in $J_\lambda^{(\alpha)}$:
$$
J_\lambda^{(\alpha)}=\sum_{\rho :|\rho|=|\lambda|} \theta^\lambda_\rho(\alpha) p_\rho.
$$
By orthogonality relations for Jack and power-sum functions,
we have its dual identity
\begin{equation} \label{eq:PowerSumJack}
p_\rho= \alpha^{\ell(\rho)} z_\rho \sum_{\lambda: |\lambda|=|\rho|} \frac{\theta^\lambda_\rho(\alpha)}{j_\lambda^{(\alpha)}} J^{(\alpha)}_\lambda
\end{equation}
and the orthogonality relation for $\theta^\lambda_\rho(\alpha)$
\begin{equation} \label{eq:OrthoTheta}
\sum_{\lambda \vdash n} \theta^\lambda_\rho(\alpha) \theta^\lambda_\pi (\alpha) \PP_n^{(\alpha)}(\lambda)=
\delta_{\rho \pi} \frac{\alpha^{n-\ell(\rho)} n!}{z_\rho}.
\end{equation}
We set $\theta^\lambda_{\mu+(1^{n-|\mu|})}(\alpha)=0$ unless $|\mu|+\ell(\mu) \le n$.
Note $\theta^{\lambda}_{(1^{|\lambda|})}(\alpha)=1$.

Let $X$ be an indeterminate.
Let $\epsilon_X$ be the algebra homomorphism from
$\mathbb{S}$ to $\mathbb{C}[X]$,
defined by
$\epsilon_X(p_r)=X$ for all $r \ge 1$.
Then we have (\cite[VI (10.25)]{Mac})
\begin{equation} \label{eq:JackSpecial}
\epsilon_X (J_\lambda^{(\alpha)}) = \prod_{(i,j) \in \lambda} (X+\alpha(j-1)-(i-1)).
\end{equation}

\subsection{Jack-Plancherel averages}

Let $A_\lambda^{(\alpha)}$ be the alphabet
$$
A_\lambda^{(\alpha)} = \{ (j-1)-(i-1)/\alpha \ | \ (i,j) \in \lambda \}.
$$
For example, $A_{(2,2)}^{(\alpha)}=\{1, 0,-1/\alpha,1-1/\alpha\}$.
Note that $A_\lambda=A^{(1)}_\lambda$ and $A_\lambda'=\{2 z \ | \ z \in A^{(2)}_\lambda\}$,
which are defined in Subsection \ref{subsec:partitions}.

Given a symmetric function $F$, we define
$$
\AV_0^{(\alpha)}(F,n)= 
\alpha^{\deg F}\sum_{\lambda \vdash n} F(A_\lambda^{(\alpha)}) \PP^{(\alpha)}_n(\lambda).
$$
More generally, for a partition $\mu$, we define 
$$
\AV_\mu^{(\alpha)}(F,n)= 
\frac{\alpha^{\deg F-|\mu|}z_{\mu+(1^{n-|\mu|})}}{ n!}\sum_{\lambda \vdash n} 
F(A_\lambda^{(\alpha)}) \PP^{(\alpha)}_n(\lambda)  \theta^\lambda_{\mu+(1^{n-|\mu|})}(\alpha).
$$
Note that $\AV_0^{(\alpha)}(F,n)= \AV_{\mu}^{(\alpha)}(F,n)$ with $\mu=(0)$.

If $F$ is homogeneous, then $F(A_\lambda^{(\alpha)})= (-\alpha)^{-\deg F} 
F(A_{\lambda'}^{(\alpha^{-1})})$.
The $\theta^\lambda_{\mu+(1^{|\lambda|-|\mu|})}(\alpha)$ has the duality
$\theta^\lambda_{\mu+(1^{|\lambda|-|\mu|})}(\alpha)= (-\alpha)^{|\mu|} 
\theta^{\lambda'}_{\mu+(1^{|\lambda|-|\mu|})}(\alpha^{-1})$ (\cite[VI (10.30)]{Mac}).
Hence we have the duality relation for $\AV_\mu^{(\alpha)}(F,n)$ with a homogeneous symmetric 
function $F$:
\begin{equation} \label{eq:dualityAValpha}
\AV_\mu^{(\alpha)}(F,n)= (-\alpha)^{\deg F-|\mu|} \AV_\mu^{(\alpha^{-1})}(F,n).
\end{equation}

The following two examples give the connection to Jucys-Murphy elements.
The average $\AV_\mu^{(\alpha)}(F,n)$ with $\alpha>0$ is
a generalization of coefficients 
$L^\lambda_\mu(n), M^\lambda_\mu(n), F^k_\mu(n)$, and $G^k_\mu(n)$,
which are studied in \cite{MN} and in the first half of the present paper. 

\begin{example}[$\alpha=1$] \label{ex:alpha1AV}
From  $J^{(1)}_\lambda=H_\lambda s_\lambda$ and from
the Frobenius formula $s_\lambda= \sum_\rho z_\rho^{-1} \chi^\lambda_\rho p_\rho$,
we have 
$\theta^\lambda_\rho(1)= z_\rho^{-1} H_{\lambda} \chi^\lambda_\rho$, 
and hence
$$
\AV_\mu^{(1)}(F,n)= \sum_{\lambda \vdash n} F(A_\lambda^{(1)}) \PP^{(1)}_n(\lambda) 
\frac{H_\lambda \chi^\lambda_{\mu+(1^{n-|\mu|})}(\alpha)}{ n!} = 
\sum_{\lambda \vdash n} F(A_\lambda) \frac{f^\lambda \chi^\lambda_{\mu+(1^{n-|\mu|})}}{n!}.
$$
In particular, we have
$\AV_0^{(1)}(F,n)=
\sum_{\lambda \vdash n} F(A_\lambda^{(1)}) \PP^{(1)}_n(\lambda)$,
the average of $F(A_\lambda)$ with respect to the Plancherel measure $\PP_n^{(1)}$.
By results in \cite{MN} (see also Subsection \ref{subsec:ClassExp} in the present paper), 
we obtain the identity
\begin{equation}
F(J_1,J_2,\dots,J_n) = \sum_\mu \AV_\mu^{(1)}(F,n) \mf{c}_\mu(n).
\end{equation}
In particular, $L^\lambda_\mu(n)=\AV_\mu^{(1)}(m_\lambda,n)$ and 
$F^k_\mu(n)=\AV_\mu^{(1)}(h_k,n)$. 
\end{example}

\begin{example}[$\alpha=2$] \label{ex:alpha2AV}
Since $J^{(2)}_\lambda=Z_\lambda$, or since $\theta^\lambda_\rho(2)= 2^{|\rho|-\ell(\rho)} |\rho|! z_{\rho}^{-1} \omega^\lambda_\rho$, we have
$$
\AV_\mu^{(2)}(F,n)=2^{\deg F}\sum_{\lambda \vdash n} F(A_\lambda^{(2)}) \PP^{(2)}_n(\lambda) 
\omega^\lambda_{\mu+(1^{n-|\mu|})} = 
\sum_{\lambda \vdash n} F(A_\lambda') \frac{f^{2\lambda}\omega^\lambda_{\mu+(1^{n-|\mu|})}}{(2n-1)!!}. 
$$
Now the equation \eqref{eq:Msum} implies 
$M^\lambda_\mu(n)=\AV_\mu^{(2)}(m_\lambda,n)$ and 
$G^k_\mu(n)=\AV_\mu^{(2)}(h_k,n)$.
More generally,
\begin{equation}
F(J_1,J_3,\dots,J_{2n-1}) \cdot P_n = \sum_\mu \AV_\mu^{(2)}(F,n) \psi_\mu(n).
\end{equation}
\end{example}

\subsection{The $\alpha=1/2$ case}

We construct the $\alpha=1/2$ version of Example \ref{ex:alpha1AV} and Example \ref{ex:alpha2AV}.
We refer to \cite[VII.2, Example 6,7]{Mac}.
Let $\epsilon$ denote the sign character of $S_{2n}$, and 
let $\epsilon_n$ denote its restriction to $H_n$:
$\epsilon_n= \epsilon \downarrow^{S_{2n}}_{H_n}$.
Then $(S_{2n},H_n, \epsilon_n)$ is a twisted Gelfand pair 
in the sense of \cite[VII.1, Example 10]{Mac}.
The corresponding Hecke algebra is
$$
\mcal{H}_n^{\epsilon}=\{ f:S_{2n} \to \mathbb{C} \ | \ 
f(\zeta \sigma)= f(\sigma \zeta) = \epsilon_n(\zeta) f(\sigma) \quad
(\sigma \in S_{2n}, \  \zeta \in H_n)\}.
$$
This algebra is commutative.
For each partition $\lambda \vdash n$,
the \emph{$\epsilon$-spherical function} $\pi^{\lambda}$ is defined by
$$
\pi^{\lambda}= (2^n n!)^{-1}\chi^{\lambda \cup \lambda} \cdot P_n^{\epsilon}
=(2^n n!)^{-1} P_n^{\epsilon} \cdot \chi^{\lambda \cup \lambda},
$$
where $P_n^{\epsilon}=\sum_{\zeta \in H_n} \epsilon_n(\zeta )\zeta$.

For each $f \in \mcal{H}_n$, let $f^{\epsilon}$ be the function on $S_{2n}$
defined by $f^{\epsilon} (\sigma)=\epsilon(\sigma) f(\sigma)$.
Then the map $f \mapsto f^{\epsilon}$ is an isomorphism of $\mcal{H}_n$ to 
$\mcal{H}_n^{\epsilon}$.
Under this isomorphism, $P_n$, $\omega^\lambda$, and $\psi_\mu(n)$ are mapped to
$P_n^{\epsilon}$, $\pi^{\lambda'}$, and $\psi_\mu^{\epsilon}(n)=
\sum_{\sigma \in H_{\mu+(1^{n-|\mu|})}} \sgn(\sigma) \sigma$, respectively.
Furthermore, for any homogeneous symmetric function $F$,
the element $F(J_1,J_3,\dots,J_{2n-1}) \cdot P_n \in \mcal{H}_n$ is mapped to
$(-1)^{\deg F} F(J_1,J_3,\dots,J_{2n-1}) \cdot P_n^{\epsilon}$,
and we have $F(A_{\lambda}')= (-1)^{\deg F} F(A_{\lambda'}^{(1/2)})$.
Therefore we can obtain the following statements from facts for 
$F(J_1,J_3,\dots,J_{2n-1}) \cdot P_n$.

\begin{thm}
\begin{enumerate}
\item For each $0 \le k< n$, we have 
$$
e_k(J_1,J_3,\dots,J_{2n-1}) \cdot P_n^{\epsilon} = (-1)^k \sum_{\mu \vdash n} 
\psi_\mu^{\epsilon}(n).
$$
\item For any symmetric function $F$ and partition $\lambda$ of $n$,
$$
F(J_1,J_3,\dots,J_{2n-1}) \cdot \pi^\lambda= \pi^\lambda \cdot 
F(J_1,J_3,\dots,J_{2n-1}) = F(A_{\lambda}^{(1/2)}) \pi^\lambda.
$$
This belongs to $\mcal{H}_n^{\epsilon}$.
\item For any symmetric function $F$,
$$
F(J_1,J_3,\dots,J_{2n-1}) \cdot P_n^{\epsilon} = \frac{1}{(2n-1)!!} \sum_{\lambda \vdash n}
f^{\lambda \cup \lambda} F(A^{(1/2)}_\lambda) \pi^\lambda.
$$
Furthermore, if $F$ is homogeneous, then
\begin{align*}
F(J_1,J_3,\dots,J_{2n-1}) \cdot P_n^{\epsilon}
=& (-1)^{\deg F} \sum_{|\mu|+\ell(\mu) \le n} \AV_\mu^{(2)}(F,n) \psi_\mu^{\epsilon}(n) \\
=&  \sum_{|\mu|+\ell(\mu) \le n} (-1)^{|\mu|}2^{\deg F-|\mu|}
\AV_\mu^{(1/2)}(F,n) \psi_\mu^{\epsilon}(n).
\end{align*}
\item For any homogeneous symmetric function $F$,
$$
\AV_\mu^{(1/2)}(F,n)= \frac{(-1)^{|\mu|}}{2^{\deg F-|\mu|}} 
\sum_{\lambda \vdash n} F(A^{(1/2)}_\lambda) 
\frac{f^{\lambda \cup \lambda} \omega^{\lambda'}_{\mu+(1^{n-|\mu|})}}{(2n-1)!!}.
$$
\end{enumerate}
\end{thm}

In a similar way to Section \ref{section:WgOrtho},
we can observe the deep connection between $F(J_1,J_3,\dots,J_{2n-1}) \cdot P_n^{\epsilon}$
and integrals over symplectic groups.
That connection will be seen in the forthcoming paper.

\subsection{Some properties}

\begin{lem} \label{lem:EvaJackContent}
Let $F$ be a symmetric function and  $n$ a positive integer. 
Assume that there exist complex numbers $\{a(\mu) \ | \ \text{$\mu$ is a partition}\}$ 
such that 
$$
F(A_\lambda^{(\alpha)})= \sum_{\mu} a(\mu) \theta^\lambda_{\mu+(1^{|\lambda|-|\mu|})}(\alpha)
$$
for any partitions $\lambda$.
Then $\AV_\mu^{(\alpha)}(F,n)= \alpha^{\deg F} a(\mu)$ for each $\mu$.
\end{lem}

\begin{proof}
We have
$$
\AV_\mu^{(\alpha)}(F,n) = 
\frac{\alpha^{\deg F- |\mu|}z_{\mu+(1^{n-|\mu|})}}{ n!} 
\sum_{\nu} a(\nu) 
\sum_{\lambda \vdash n} 
\PP^{(\alpha)}_n(\lambda)  \theta^\lambda_{\mu+(1^{n-|\mu|})}(\alpha)
\theta^\lambda_{\nu+(1^{n-|\nu|})}(\alpha).
$$
The claim follows from the orthogonality relation \eqref{eq:OrthoTheta}.
\end{proof}

The following theorem is a Jack deformation of
Jucys' result \cite{Jucys} and its analogue, Corollary \ref{cor1-Jucys}.

\begin{prop} \label{prop:JucysAlpha}
$$
\AV^{(\alpha)}_\mu(e_k,n) 
=\begin{cases} 1 & \text{if $|\mu|=k$}, \\
0 & \text{otherwise}.
\end{cases}
$$
\end{prop}

\begin{proof}
Let $X$ be an indeterminate and let $\lambda \vdash n$.
It follows from \eqref{eq:JackSpecial}  that
\begin{align*}
& \sum_{k=0}^n e_k(A_\lambda^{(\alpha)}) X^k =
(X/\alpha)^n \prod_{(i,j) \in \lambda}  (\alpha /X +\alpha(j-1)-(i-1)) 
= (X/\alpha)^n \epsilon_{\alpha/X} (J_\lambda^{(\alpha)}) \\
=& (X/\alpha)^n \sum_{\rho \vdash n} \theta^\lambda_\rho(\alpha) 
\epsilon_{\alpha/X} (p_\rho)  
= \sum_{\rho \vdash n} \theta^\lambda_\rho(\alpha)  (X/\alpha)^{n-\ell(\rho)}
=  \sum_{k=0}^n \alpha^{-k} 
\sum_{\nu \vdash k} \theta^\lambda_{\nu+(1^{n-k})}(\alpha)  X^{k},
\end{align*}
which gives
\begin{equation} \label{eq:La98}
e_k(A_\lambda^{(\alpha)}) =\alpha^{-k} \sum_{\nu \vdash k} 
\theta^\lambda_{\nu+(1^{|\lambda|-k})}(\alpha).
\end{equation}
The claim follows from this identity and Lemma \ref{lem:EvaJackContent}.
\end{proof}

\begin{remark}
The equation \eqref{eq:La98} is seen in \cite[Theorem 5.4]{Lassalle98}.
\end{remark}

\begin{thm} \label{thm:AVpolynomial}
Let $F$ be any symmetric function  and $\mu$ a partition.
Then   
$\AV_\mu^{(\alpha)}(F,n)$ is a polynomial in $n$.
If the expansion of $F$ in $p_\rho$ is given by
$F=\sum_\rho a(\rho) p_\rho$, then
the degree of $\AV_\mu^{(\alpha)}(F,n)$ in $n$ is at most
$$
\max_{a(\rho) \not=0} (|\rho|+\ell(\rho)) - (|\mu|+\ell(\mu)).
$$
\end{thm}

This theorem at $\alpha=2$ with Example \ref{ex:alpha2AV} implies 
part 1 of Theorem \ref{thm:coefM}.
The proof is given in the next subsection by applying shifted symmetric function theory.

\begin{example}
Since the  monomial symmetric function is expanded as 
$$
m_\lambda= p_\lambda +\sum_{\rho > \lambda} a(\rho) p_\rho, 
$$
the degree of the polynomial $\AV^{(\alpha)}_\mu(m_\lambda,n)$ is 
at most $|\lambda|+\ell(\lambda) -(|\mu|+\ell(\mu))$.
But this evaluation is not sharp. 
Indeed, as we will observe below, 
the degree of $\AV_{(0)}^{(\alpha)}(m_{(3)},n)=\alpha(\alpha-1) \binom{n}{2}$ is $2$
but $(|\lambda|+\ell(\lambda))-(|\mu|+\ell(\mu))=4$ with $\lambda=(3)$ and $\mu=(0)$.
\end{example}

\subsection{Shifted symmetric functions and proof of Theorem \ref{thm:AVpolynomial}}

Following to \cite{KOO, LassalleSomeIdentities},
we review the theory of shifted symmetric functions related to 
Jack functions.

A polynomial  in $n$ variables $x_1,x_2,\dots,x_n$ is said to be 
\emph{shifted-symmetric} if it is symmetric in the variables
$y_i:=x_i-i/\alpha$.
Denote by $\mathbb{S}^*_\alpha(n)$ the subalgebra of
shifted-symmetric functions in $\mathbb{C}[x_1,x_2,\dots,x_n]$.

Consider an infinite alphabet $x=(x_1,x_2,\dots)$
and consider the morphism $F(x_1,x_2,\dots,x_n,x_{n+1}) \mapsto
F(x_1,x_2,\dots,x_{n},0)$ from $\mathbb{S}^*_\alpha(n+1)$ to $\mathbb{S}^*_\alpha(n)$.
As  the definition of $\mathbb{S}$,
we can define the algebra $\mathbb{S}^*_\alpha$ as the 
projective limit of the sequence $(\mathbb{S}_\alpha^*(n))_{n \ge 1}$.
Elements of $\mathbb{S}_\alpha^*$ are called \emph{shifted-symmetric functions}
and written as $F(x)=F(x_1,x_2,\dots)$ using infinite variables.
Denote by $\deg F$ the degree of $F$.

For each $F \in \mathbb{S}^*_\alpha$, we may evaluate at
partitions:
$F(\lambda)=F(\lambda_1,\lambda_2,\dots)$.
We denote by $[F] \in \mathbb{S}$ the homogeneous symmetric terms of degree $\deg F$.
We call $[F]$ the \emph{leading symmetric term} of $F$.
The map $F \mapsto [F]$ provides a canonical isomorphism of the graded algebra associated 
to the filtered algebra $\mathbb{S}^*_\alpha$ onto $\mathbb{S}$.
Assuming that the leading terms $[F_1], [F_2],\dots$ of a sequence $F_1,F_2,\dots$ 
in $\mathbb{S}^*_\alpha$ generate the algebra $\mathbb{S}$,
this sequence itself generates $\mathbb{S}^*_\alpha$.

For each integer $k \ge 1$, consider a polynomial
$$
p^*_k(x;\alpha)= \sum_{i \ge 1} 
\( (x_i-(i-1)/\alpha)^{\downarrow k}- (-(i-1)/\alpha)^{\downarrow k} \)
$$
with $a^{\downarrow k}=a(a-1) \cdots (a-k+1)$.
Then these polynomials are shifted-symmetric.
Since $[p^*_k(\cdot;\alpha)]=p_k$
and since the $p_k$ generate $\mathbb{S}$, 
they generate $\mathbb{S}^*_\alpha$.

For $F \in \mathbb{S}$ and a partition $\lambda$,
we put $H_F^{(\alpha)}(\lambda)=F(A_\lambda^{(\alpha)})$.

\begin{lem}[Lemma 7.1 in \cite{LassalleSomeIdentities}]
For any integer $k \ge 1$, the function $\lambda \mapsto H_{p_k}^{(\alpha)}(\lambda)
=p_k(A^{(\alpha)}_\lambda)$
defines a shifted symmetric function of degree $\deg H_{p_k}^{(\alpha)}=k+1$.
Specifically,
$$
H_{p_k}^{(\alpha)}(\lambda)= \sum_{m=1}^k S(k,m) \frac{p^*_{m+1}(\lambda;\alpha)}{m+1}.
$$
Here $S(k,m)$ are Stirling's numbers of second kinds,
defined via $u^k=\sum_{m=1}^k S(k,m) u^{\downarrow m}$.
\end{lem}

Since $p_k$ generate $\mathbb{S}$, we have the following corollary.

\begin{cor} \label{cor:DegreeHF}
For any $F \in \mathbb{S}$, the function $\lambda \mapsto H_F^{(\alpha)}(\lambda)$
defines a shifted symmetric function.
Furthermore, if the expansion of $F$ in $p_\rho$ is given by
$F=\sum_\rho a(\rho) p_\rho$, then the degree of $H_F^{(\alpha)}$ is 
$\max_{a(\rho)\not=0} (|\rho|+\ell(\rho))$.
\end{cor}

Now we define \emph{shifted Jack functions} $J_\mu^*(x;\alpha)$.
They are defined by
$$
J^*_\mu(\lambda;\alpha)= \frac{|\lambda|^{\downarrow |\mu|} 
\langle p_1^{|\lambda|-|\mu|} J_\mu^{(\alpha)}, J_\lambda^{(\alpha)} \rangle_\alpha}
{\alpha^{|\lambda|} |\lambda|!}.
$$

\begin{lem}[\cite{KOO}] \label{lem:PropertySJack}
The $J^*_\mu(x;\alpha)$ are shifted-symmetric and satisfy the following properties.
\begin{enumerate}
\item $[J^*_\mu(\cdot;\alpha)]=J^{(\alpha)}_\mu$.
Hence the $J^*_\mu(x;\alpha)$ form a basis of $\mathbb{S}^*_\alpha$.
\item $J_\mu^*(\lambda;\alpha)=0$ unless $\mu_i \le \lambda_i$ for all $i \ge 1$.
\item $J_\mu^*(\mu;\alpha)= \alpha^{-|\mu|} j_\mu^{(\alpha)}$.
\end{enumerate}
\end{lem}

The following theorem is a slight extension of Theorem 5.5 in \cite{Olshanski}.

\begin{prop} \label{prop:AVsJack}
Let $\mu,\nu$ be partitions.
If $|\nu| \ge |\mu|+\ell(\mu)$, then we have 
%\begin{equation} \label{eq:GeneralOlshanski}
$$
\frac{z_{\mu+(1^{n-|\mu|})}}{ n!}\sum_{\lambda \vdash n} 
J_\nu^*(\lambda;\alpha) \PP^{(\alpha)}_n(\lambda)  \theta^\lambda_{\mu+(1^{n-|\mu|})}(\alpha)=
\binom{n-|\mu|-\ell(\mu)}{|\nu|-|\mu|-\ell(\mu)} z_{\mu+(1^{|\nu|-|\mu|})}
\theta^\nu_{\mu+(1^{|\nu|-|\mu|})}(\alpha),
$$
which is a polynomial in $n$ of degree $|\nu|-|\mu|-\ell(\mu)$.
Otherwise, both sides are zero.
\end{prop}

\begin{proof}
Put $m=|\nu|$.
If $n <m$, then both sides vanish by part 2 of Lemma \ref{lem:PropertySJack},
and so we may assume $n \ge m$. 
We have
\begin{align*}
& \frac{z_{\mu+(1^{n-|\mu|})}}{ n!}\sum_{\lambda \vdash n} 
J_\nu^*(\lambda;\alpha) \PP^{(\alpha)}_n(\lambda)  \theta^\lambda_{\mu+(1^{n-|\mu|})}(\alpha) \\
=& \frac{z_{\mu+(1^{n-|\mu|})}}{ n!}\sum_{\lambda \vdash n} 
\frac{n^{\downarrow m} 
\langle p_1^{n-m} J_\nu^{(\alpha)}, J_\lambda^{(\alpha)} \rangle_\alpha}
{j_\lambda^{(\alpha)} }
  \theta^\lambda_{\mu+(1^{n-|\mu|})}(\alpha) \\
 =& \frac{n^{\downarrow m}}{\alpha^{n-|\mu|} n!} 
 \left\langle p_1^{n-m} J_\nu^{(\alpha)},  
 \sum_{\lambda \vdash n} 
 \frac{\alpha^{n-|\mu|} z_{\mu+(1^{n-|\mu|})}}{j_\lambda^{(\alpha)}}
 \theta^\lambda_{\mu+(1^{n-|\mu|})} J_\lambda^{(\alpha)} \right\rangle_\alpha \\
 =& \frac{n^{\downarrow m}}{\alpha^{n-|\mu|} n!} 
 \left\langle p_1^{n-m} J_\nu^{(\alpha)},  p_{\mu+(1^{n-|\mu|})} \right\rangle_\alpha
 \qquad \text{by \eqref{eq:PowerSumJack}}.
\end{align*}
Using the fact that the adjoint to the multiplication by $p_1$ 
with respect to $\langle \cdot, \cdot \rangle_\alpha$ is 
$\alpha \frac{\partial}{\partial p_1}$, we have 
\begin{equation} \label{eq:Polynomiality1}
=\frac{n^{\downarrow m}}{\alpha^{n-|\mu|} n!} 
 \left\langle  J_\nu^{(\alpha)}, \alpha^{n-m} \( \frac{\partial}{\partial p_1}\)^{n-m}
  p_{\mu+(1^{n-|\mu|})} \right\rangle_\alpha.
\end{equation}
Since $m_1(\mu+(1^{n-|\mu|}))=n-|\mu|-\ell(\mu)$, 
the symmetric function $\(\tfrac{\partial}{\partial p_1} \)^{n-m} p_{\mu+(1^{n-|\mu|})}$ 
vanishes unless $m \ge |\mu|+\ell(\mu)$.
If $n \ge m \ge |\mu|+\ell(\mu)$, then \eqref{eq:Polynomiality1} equals
\begin{align*}
=& \frac{n^{\downarrow m}}{\alpha^{m-|\mu|} n!} (n-|\mu|-\ell(\mu))^{\downarrow (n-m)} 
\langle J_\nu^{(\alpha)}, p_{\mu+(1^{m-|\mu|})} \rangle_\alpha \\
=& \binom{n-|\mu|-\ell(\mu)}{m-|\mu|-\ell(\mu)} z_{\mu+(1^{m-|\mu|})}
\theta^\nu_{\mu+(1^{m-|\mu|})}(\alpha).
\end{align*}
\end{proof}

\begin{remark}
Proposition \ref{prop:AVsJack} can be rewritten as follows: 
for partitions $\nu,\mu$ such that $|\nu| \ge |\mu|+\ell(\mu)$ and for any $n \ge 0$,
$$
\sum_{\lambda \vdash n} 
J_\nu^*(\lambda;\alpha) \PP^{(\alpha)}_n(\lambda)  \theta^\lambda_{\mu+(1^{n-|\mu|})}(\alpha)=
n^{\downarrow |\nu|}\theta^\nu_{\mu+(1^{|\nu|-|\mu|})}(\alpha).
$$
In particular, we obtain a simple identity
$$
\sum_{\lambda \vdash n} 
J_\nu^*(\lambda;\alpha) \PP^{(\alpha)}_n(\lambda)=
n^{\downarrow |\nu|},
$$
which is seen in \cite[Theorem 5.5]{Olshanski}.
Lassalle obtained a similar identity.
Specifically, Equation (3.3) in \cite{Lassalle08} implies that,
for partitions $\nu$ and $\mu$ such that $|\nu| \ge |\mu|+\ell(\mu)$
and for any $m \ge 0$, 
$$
 \sum_{\rho \vdash m} J^*_\rho(\nu;\alpha)
\PP_m^{(\alpha)}(\rho) \theta^\rho_{\mu +(1^{m-|\mu|})}(\alpha)=
\frac{(|\nu|-|\mu|-\ell(\mu))! \, m!}{(|\nu|-m)! \, (m-|\mu|-\ell(\mu))!} 
\theta^\nu_{\mu+(1^{|\nu|-|\mu|})}(\alpha).
$$
\end{remark}

\begin{proof}[Proof of Theorem \ref{thm:AVpolynomial}
and part 1 of Theorem \ref{thm:coefM}]
The statement follows from Theorem \ref{cor:DegreeHF}, part 1 of Lemma \ref{lem:PropertySJack},
and Proposition \ref{prop:AVsJack}.
\end{proof}

\section{Examples and open problems} 

\subsection{Examples of $\AV_\mu^{(\alpha)}(F,n)$}
\label{subsec:Examples}

We give examples of $\AV_\mu^{(\alpha)}(m_\lambda,n)$ and
$\AV_\mu^{(\alpha)}(h_k,n)$, studied in the previous section.

\noindent
$|\lambda|=0,1$.
$$
\AV_{\mu}^{(\alpha)}(m_{(0)},n)= \delta_{\mu,(0)}. \qquad
\AV_{\mu}^{(\alpha)}(m_{(1)},n)= \delta_{\mu,(1)}.
$$
\noindent
$|\lambda|=2$.
\begin{align*}
\AV_\mu^{(\alpha)}(m_{(2)},n)=&\delta_{\mu,(2)}+(\alpha-1) \delta_{\mu,(1)} +\alpha \binom{n}{2}
\delta_{\mu,(0)}. \\
\AV_\mu^{(\alpha)}(m_{(1^2)},n)=&\delta_{\mu,(2)}+\delta_{\mu,(1^2)}. \\
\AV_\mu^{(\alpha)}(h_2,n)=&2\delta_{\mu,(2)}+\delta_{\mu,(1^2)}
+(\alpha-1) \delta_{\mu,(1)} +\alpha \binom{n}{2}
\delta_{\mu,(0)}.
\end{align*}
$|\lambda|=3$.
\begin{align*}
\AV_\mu^{(\alpha)}(m_{(3)},n)=&
\delta_{\mu,(3)}+
3(\alpha-1)\delta_{\mu,(2)}+
(2\alpha n +\alpha^2-5\alpha+1)\delta_{\mu,(1)} 
+\alpha(\alpha-1) \binom{n}{2}\delta_{\mu,(0)}. \\
\AV_\mu^{(\alpha)}(m_{(2,1)},n)=&
3\delta_{\mu,(3)}+ \delta_{\mu,(2,1)} 
+3(\alpha-1)\delta_{\mu,(2)}+ 2(\alpha-1) \delta_{\mu,(1^2)}
+\alpha \(\binom{n}{2}-1\)\delta_{\mu,(1)}. \\
\AV_\mu^{(\alpha)}(m_{(1^3)},n)=&
\delta_{\mu,(3)}+\delta_{\mu,(2,1)}+\delta_{\mu,(1^3)}. \\
\AV_\mu^{(\alpha)}(h_3,n)=&
5\delta_{\mu,(3)}+2\delta_{\mu,(2,1)}+\delta_{\mu,(1^3)}+
6(\alpha-1)\delta_{\mu,(2)}+2(\alpha-1)\delta_{\mu,(1^2)} \\
&\qquad +
\( \frac{1}{2} \alpha n^2 + \frac{3}{2} \alpha n+\alpha^2-6\alpha+1\)\delta_{\mu,(1)} 
+\alpha(\alpha-1) \binom{n}{2}\delta_{\mu,(0)}. 
\end{align*}

In fact,
the identities for $m_{(1^k)}$ are given by Proposition \ref{prop:JucysAlpha}.
Lassalle \cite{LassalleCumulant, LassalleHP} (see also \cite[Conjecture 8.1]{Lassalle98}) gives the expansion of 
$\theta^\lambda_{\mu+(1^{|\lambda|-|\mu|})}(\alpha)$ with respect 
to $p_\rho(A^{(\alpha)}_\lambda)$:
Letting  
$\hat{\theta}_\mu(\lambda)=\theta^\lambda_{\mu+(1^{|\lambda|-|\mu|})}(\alpha)$ and 
$\hat{p}_\rho(\lambda)=p_\rho(A^{(\alpha)}_\lambda)$,
\begin{align*}
\hat{\theta}_{(1)}(\lambda)=& \alpha \hat{p}_{(1)}(\lambda), \\
\hat{\theta}_{(2)}(\lambda)=& \alpha^2 \hat{p}_{(2)}(\lambda) 
-\alpha(\alpha-1) \hat{p}_{(1)} (\lambda) - \alpha \binom{|\lambda|}{2},\\
\hat{\theta}_{(1^2)}(\lambda)=& -\tfrac{3}{2} \alpha^2 \hat{p}_{(2)}(\lambda) 
+\tfrac{1}{2} \alpha^2 \hat{p}_{(1^2)} (\lambda)
+\alpha(\alpha-1) \hat{p}_{(1)}(\lambda) +\alpha \binom{|\lambda|}{2}, \\
\hat{\theta}_{(3)}(\lambda)=& 
\alpha^3 \hat{p}_{(3)}(\lambda) -3\alpha^2 (\alpha-1) \hat{p}_{(2)}(\lambda) +
\alpha (-2\alpha |\lambda|+2\alpha^2-\alpha+2) \hat{p}_{(1)}(\lambda)
+2\alpha(\alpha-1) \binom{|\lambda|}{2}, \\
\hat{\theta}_{(2,1)}(\lambda)=& 
-4 \alpha^3 \hat{p}_{(3)}(\lambda) +\alpha^3 \hat{p}_{(2,1)}(\lambda) 
+9\alpha^2 (\alpha-1) \hat{p}_{(2)} (\lambda)-\alpha^2 (\alpha-1) \hat{p}_{(1^2)}(\lambda) \\
& +\alpha(-\tfrac{\alpha}{2} |\lambda|^2 +\tfrac{13 \alpha}{2}|\lambda|-5\alpha^2+2\alpha-5) \hat{p}_{(1)}(\lambda)
-5\alpha(\alpha-1) \binom{|\lambda|}{2}.
\end{align*}
Using these, we can express $p_\rho(A^{(\alpha)}_\lambda)$
in terms of $\hat{\theta}_\mu(\lambda)$ and therefore also 
$m_\nu(A^{(\alpha)}_\lambda)$.
Hence our above examples follow by Lemma \ref{lem:EvaJackContent}.

Those examples are reduced as typical cases
\begin{align*}
L^\lambda_\mu(n)=&\AV_\mu^{(1)}(m_\lambda,n), & 
F^k_\mu(n)=&\AV_\mu^{(1)}(h_k,n), \\
M^\lambda_\mu(n)=&\AV_\mu^{(2)}(m_\lambda,n), &
G^k_\mu(n)=&\AV_\mu^{(2)}(h_k,n).
\end{align*}
For example, we obtain
$$
\AV_\mu^{(2)}(h_3,n)=
5\delta_{\mu,(3)}+2\delta_{\mu,(2,1)}+\delta_{\mu,(1^3)}+
6\delta_{\mu,(2)}+2\delta_{\mu,(1^2)} 
(n^2 +3 n-7)\delta_{\mu,(1)} 
+n(n-1)\delta_{\mu,(0)},
$$
or, by Example \ref{ex:alpha2AV},
\begin{align*}
h_{3}(J_1,J_3,\dots,J_{2n-1}) \cdot P_n =& 5\psi_{(3)}(n)+2\psi_{(2,1)}(n)+\psi_{(1^3)}(n)+ 6\psi_{(2)}(n) +2 \psi_{(1^2)}(n)& \\
&\qquad +(n^2+3n-7)\psi_{(1)}(n)+n(n-1) \psi_{(0)}(n).
\end{align*}
See also the $\alpha=1$ cases in \cite[\S 8.1]{MN}.

We remark that 
a conjecture for $\AV_0^{(\alpha)}(p_\lambda,n)$ is given by
\cite[Conjecture 12.1]{Lassalle98}.

\subsection{Table of asymptotic expansions of $\mr{Wg}^{O(N)}_{n}$} \label{subsec:tableAsymWg}

We give some examples of the expansion \eqref{eq:AsymExpWg},
by using coefficients appearing in Subsection \ref{subsec:Examples}.

Given a partition $\mu$,
we define 
$\mr{Wg}^{O(N)}(\mu;n)=\mr{Wg}^{O(N)}_n(\sigma)$,
where $\sigma$ is a permutation in $S_{2n}$ of reduced coset-type $\mu$.
For example,
$$
\mr{Wg}^{O(N)}((0);n)=\mr{Wg}^{O(N)}_n(\mr{id}_{2n}), \qquad
\mr{Wg}^{O(N)}((1);n)=\mr{Wg}^{O(N)}_n\( \(\begin{smallmatrix} 
1 & 2 & 3 & 4 & 5 & 6 & \dots & 2n-1 & 2n \\
1 & 4 & 2 & 3 & 5 & 6 & \dots & 2n-1 & 2n
\end{smallmatrix} \)\).
$$

Theorem \ref{thm:ExpansionWg} and examples in the previous subsection  with $\alpha=2$
give the following 
asymptotic expansions.
As $N \to \infty$,
\begin{align*}
\mr{Wg}^{O(N)}((0);n)=& 
N^{-n} +n(n-1)N^{-n-2} -n(n-1)N^{-n-3} + \mr{O}(N^{-n-4}). \\
\mr{Wg}^{O(N)}((1);n)=& 
-N^{-n-1} +N^{-n-2} -(n^2+3n-7)N^{-n-3} + \mr{O}(N^{-n-4}). \\
\mr{Wg}^{O(N)}((2);n)=&
2N^{-n-2}-6N^{-n-3}+ \mr{O}(N^{-n-4}).\\
\mr{Wg}^{O(N)}((1^2);n)=&
N^{-n-2}-2 N^{-n-3}+ \mr{O}(N^{-n-4}).
\end{align*}

On the other hand, in \cite{CM},
the explicit values of $\mr{Wg}^{O(N)}_{n}(\sigma)$ for  $n \le 6$ are given.
We remark that in \cite{CM} ordinary coset-types were used, not reduced ones.
Using a computer with the table in \cite{CM}, we obtain the following expansions.

\begin{align*}
\mr{Wg}^{O(N)}((0);2)=& N^{-2} - 0 N^{-3}+2 N^{-4}-2 N^{-5}+6 N^{-6} -10 N^{-7} +22 N^{-8} -\cdots. \\
\mr{Wg}^{O(N)}((0);3)=& N^{-3}- 0 N^{-4} +6 N^{-5}-6 N^{-6}+50 N^{-7} -126 N^{-8} +610 N^{-9} -\cdots.\\
\mr{Wg}^{O(N)}((0);4)=& N^{-4}- 0 N^{-5} +12 N^{-6}-12 N^{-7}+176 N^{-8} -468 N^{-9} +3544 N^{-10} 
-\cdots.\\
\mr{Wg}^{O(N)}((0);5)=& N^{-5} - 0 N^{-6}+20 N^{-7}-20 N^{-8}+440 N^{-9} -1180 N^{-10} +12480 N^{-11}
 -\cdots. \\
\mr{Wg}^{O(N)}((0);6)=& N^{-6}- 0 N^{-7} +30 N^{-8}-30 N^{-9}+910 N^{-10} -2430 N^{-11} +33710 N^{-12}
 -\cdots.
\end{align*}
\begin{align*}
\mr{Wg}^{O(N)}((1);2)=& -N^{-3} +N^{-4} -3 N^{-5}+5N^{-6} -11 N^{-7} +21N^{-8}-43N^{-9}+ \cdots.\\
\mr{Wg}^{O(N)}((1);3)=& -N^{-4} +N^{-5} -11 N^{-6}+29N^{-7} -147 N^{-8} +525 N^{-9}
-2227N^{-10}+ \cdots.  \\
\mr{Wg}^{O(N)}((1);4)=& -N^{-5} +N^{-6} -21 N^{-7}+57N^{-8} -489 N^{-9} +2157 N^{-10}
-14077N^{-11}+ \cdots.  \\
\mr{Wg}^{O(N)}((1);5)=& -N^{-6} +N^{-7} -33 N^{-8}+89N^{-9} -1117 N^{-10} +5237 N^{-11}
-45881N^{-12}+ \cdots.\\
\mr{Wg}^{O(N)}((1);6)=& -N^{-7} +N^{-8} -47 N^{-9}+125N^{-10} -2123 N^{-11} +10121 N^{-12}
-112551N^{-13}+ \cdots.
\end{align*}
\begin{align*}
\mr{Wg}^{O(N)}((2);3)=& 2N^{-5} -6N^{-6} +34 N^{-7}-126N^{-8} +546 N^{-9} -2142N^{-10}+ \cdots. \\
\mr{Wg}^{O(N)}((2);4)=& 2N^{-6} -6N^{-7} +64 N^{-8}-300N^{-9} +2094 N^{-10} -11682N^{-11}+ \cdots.\\
\mr{Wg}^{O(N)}((2);5)=& 2N^{-7} -6N^{-8} +98 N^{-9}-490N^{-10} +4694 N^{-11} -30382N^{-12}+ \cdots.\\
\mr{Wg}^{O(N)}((2);6)=& 2N^{-8} -6N^{-9} +136 N^{-10}-696N^{-11} +8590 N^{-12} -59850N^{-13}+ \cdots.
\end{align*}
\begin{align*}
\mr{Wg}^{O(N)}((1^2);4)=& N^{-6} -2N^{-7} +43 N^{-8}-216N^{-9} +1737 N^{-10} -10254N^{-11}+ \cdots. \\
\mr{Wg}^{O(N)}((1^2);5)=& N^{-7} -2N^{-8} +59 N^{-9}-280N^{-10} +3257 N^{-11} -21934N^{-12}+ \cdots. \\
\mr{Wg}^{O(N)}((1^2);6)=& N^{-8} -2N^{-9} +77 N^{-10}-350N^{-11} +5385 N^{-12} -37498N^{-13}+ \cdots. 
\end{align*}
\begin{align*}
\mr{Wg}^{O(N)}((3);4)=& -5N^{-7}+29N^{-8} -258 N^{-9} +1590 N^{-10} -10695 N^{-11} + \cdots. \\
\mr{Wg}^{O(N)}((3);5)=& -5N^{-8}+29N^{-9} -370 N^{-10} +2630 N^{-11} -23815 N^{-12} + \cdots.\\
\mr{Wg}^{O(N)}((3);6)=& -5N^{-9}+29N^{-10}-492 N^{-11} +3738 N^{-12} -42019 N^{-13} + \cdots.
\end{align*}
\begin{align*}
\mr{Wg}^{O(N)}((2,1);5)=& -2N^{-8}+8N^{-9}-190N^{-10}+1460 N^{-11} - 15994 N^{-12}+ \cdots. \\
\mr{Wg}^{O(N)}((2,1);6)=& -2N^{-9}+8N^{-10}-236N^{-11}+1760 N^{-12} - 24254 N^{-13}+ \cdots.
\end{align*}
$$
\mr{Wg}^{O(N)}((1^3);6)= -N^{-9}+3N^{-10}-120N^{-11}+742N^{-12}-13023N^{-13}+\cdots. 
$$
\begin{align*}
\mr{Wg}^{O(N)}((4);5)=& 14N^{-9}-130N^{-10}+1640 N^{-11} - 14740 N^{-12}+138578N^{-13}-\cdots. \\
\mr{Wg}^{O(N)}((4);6)=& 14N^{-10}-130N^{-11}+2060 N^{-12} - 20360 N^{-13}+232838N^{-14}-\cdots.
\end{align*}
$$
\mr{Wg}^{O(N)}((3,1);6)= 5N^{-10}-34N^{-11}+862N^{-12}-9096N^{-13}+126523N^{-14}+\cdots. 
$$
$$
\mr{Wg}^{O(N)}((2,2);6)= 4N^{-10}-24N^{-11}+772N^{-12}-8436N^{-13}+121936N^{-14}+\cdots. 
$$
$$
\mr{Wg}^{O(N)}((5);6)= -42N^{-11}+562N^{-12}-9426N^{-13}+114478N^{-14}-\cdots. 
$$

\subsection{Open questions} \label{subsec:Open}

\begin{enumerate}
\item (cf. Corollary \ref{cor:SymPolyHecke}.)
It is known that the set $\{F(J_1,J_2,\dots,J_n) \ | \ F \in \mbb{S} \}$
coincides with the center $\mcal{Z}_n$ of the group algebra $\bC[S_n]$.
Thus symmetric functions in Jucys-Murphy elements generate $\mcal{Z}_n$.
Now the following conjecture is natural.

\begin{conj}
The set $\{F(J_1,J_3,\dots,J_{2n-1}) \cdot P_n\ | \ F \in 
\mbb{S} \}$ coincides with the Hecke algebra $\mcal{H}_n$.
\end{conj}

\item (cf. part 4 of Theorem \ref{thm:coefM}
and Examples in  Subsection \ref{subsec:Examples}.)
We suggest the following conjecture.
\begin{conj} 
Let $F$ be a symmetric function of degree $k$ and
let $\alpha$ be  a positive real number.
Then, for each partition $\mu \vdash k$, 
$\AV^{(\alpha)}_\mu(F,n)$ is independent of both $\alpha$ and $n$.
In particular, for $\lambda,\mu \vdash k$,
$\AV^{(\alpha)}_\mu(m_\lambda,n)=L^\lambda_\mu$ 
and $\AV^{(\alpha)}_\mu(h_k,n)= \prod_{i=1}^{\ell(\mu)} \mr{Cat}_{\mu_i}$.
\end{conj}
\item (cf. Examples in  Subsection \ref{subsec:Examples}.)
We suggest the following conjecture.
\begin{conj}\label{conj:SecondOrder}
Let $F$ be a homogeneous symmetric function of degree $k$ and
let $\alpha$ be  a positive real number.
Then, for each partition $\mu \vdash k-1$, 
$\AV^{(\alpha)}_\mu(F,n)$ is independent of $n$
(but depends on $\alpha$).
\end{conj}
\item (cf. Theorem \ref{thm:ExpansionWg} and Subsection \ref{subsec:tableAsymWg}.)
Conjecture \ref{conj:SecondOrder} implies that $G^{|\mu|+1}_\mu(n)=\AV^{(2)}_\mu(h_{|\mu|+1},n)$
is independent of $n$.
Can you evaluate the $G^{|\mu|+1}_\mu=G^{|\mu|+1}_\mu(n)$ explicitly?
From identities in Subsection \ref{subsec:tableAsymWg},
we can obtain  
$$
G^{1}_{(0)} =0, \qquad G^{2}_{(1)}=1, \qquad
G^{3}_{(2)} =6, \qquad G^{3}_{(1^2)} =2,
$$
and conjecture
\begin{gather*}
G^{4}_{(3)}\stackrel{?}{=} 29, \qquad 
G^{4}_{(2,1)}\stackrel{?}{=} 8,\qquad
G^{4}_{(1^3)}\stackrel{?}{=} 3,\\
G^{5}_{(4)}\stackrel{?}{=} 130, \qquad 
G^{5}_{(3,1)}\stackrel{?}{=} 34, \qquad
G^{5}_{(2^2)}\stackrel{?}{=} 24, \qquad
G^6_{(5)} \stackrel{?}{=} 562.
\end{gather*}
Recall that the $n$-independent number $G^{|\mu|}_\mu=G^{|\mu|}_\mu(n)$ 
is the product of Catalan numbers.
How about $G^{|\mu|+1}_\mu$?
We could expect that $G^{|\mu|+1}_\mu$ has a good combinatorial interpretation.
For one-row partitions, we  suggest the following conjecture.

\begin{conj}
Let $n,k$ be nonnegative integers such that $n>k$.
Then $G^{k+1}_{(k)}(n)$ is independent of $n$ and equal to
$$
4^k - \binom{2k+1}{k}.
$$
Equivalently,
$$
\mr{Wg}^{O(N)}((k);n) \stackrel{?}{=}
(-1)^k \mr{Cat}_k N^{-n-k} +(-1)^{k+1} \(4^k-\binom{2k+1}{k}\) N^{-n-k-1} +\mr{O}(N^{-n-k-2}).
$$
\end{conj}

The number $4^k - \binom{2k+1}{k}$ is called 
{\it the area of Catalan paths} of length $k$,
see \cite{CEF}.
Define the set $\mf{E}(k)$ by
$$
\mf{E}(k)=\left\{(i_1,i_2,\dots,i_k) \in \mathbb{Z}^k \ \Bigm| \ 
\begin{array}{l}
1 \le i_1 \le i_2 \le \cdots \le i_k \le k, \\
i_p \ge p  \ (1 \le  p \le k)
\end{array}
\right\}.
$$
It is known that
$$
\mr{Cat}_k = |\mf{E}(k)|, \qquad 4^k-\binom{2k+1}{4}
= \sum_{(i_1,i_2,\dots,i_k) \in \mf{E}(k)}
 \sum_{p=1}^k (2(i_p-p)+1).
$$
\end{enumerate}

\medskip

\section*{Acknowledgements}

The author would like to acknowledge a lot of kind comments with Michel Lassalle
and thanks for reviewers' suggestions under a revision.

%<<<<<<<<<<<<<<<<<<<<<<<<<<<<<<<<<<<<<<

\end{document}